\documentclass{amsbook}
\usepackage{amssymb}
\usepackage{amsfonts}

\makeatletter
\numberwithin{equation}{section}
\numberwithin{figure}{section}
\numberwithin{equation}{section}
\numberwithin{figure}{section}
\numberwithin{equation}{section}

\theoremstyle{definition}

\newtheorem*{acknowledgements*}{Acknowledgements}
\theoremstyle{remark}

\numberwithin{theorem}{section}   
\input{tcilatex}

\begin{document}
\frontmatter
\title[Banach Couples]{Banach Couples. I. Elementary theory.}
\author{Jaak Peetre}
\address{Lund Institute of Technology, Lund, Sweden}
\author{Per G. Nilsson (Typist)}
\address{Nilsson: Stockholm, Sweden}
\email{pgn@plntx.com}
\date{\today }

\begin{abstract}
This note is an (exact) copy of the report of Jaak Peetre, "Banach Couples.
I. Elementary Theory". Published as Technical Report, Lund (1971). Some more
recent general references have been added and some references updated though.
\end{abstract}

\maketitle

\begin{center}
{\Huge Banach Couples, I. Elementary Theory.}

{\Huge Jaak Peetre}
\end{center}

\section*{0. Introduction}

The theory of interpolation spaces deals roughly speaking with
"constructions", so-called "interpolation methods", which to a given set of
Banach spaces (usually consisting of precisely two elements; we speak then
of \underline{\emph{Banach couples}}) assign a Banach space, called 
\underline{\emph{interpolation}}\emph{\ }\underline{\emph{space }}with
respect to the set of Banach spaces. When trying to put this very vague
concepts on a more rigorous basis one is lead to consider functors from a
category of sets of Banach spaces into a category of Banach spaces. I have
long felt that in order to push any further it is necessary to consider
categories of the first named type \underline{\emph{per se}}.$\medskip $%
\footnote{%
For references, post 1971, see notably Brudnyi-Krugljak-Aizenstein \cite[%
Chapter 2, Subsection 2.7.2]{49-BK}, Brudnyi-Krein-Semenov \cite{50-BrKrSe},
Kaijer-Pelletier \cite{52-KaPe}, \cite{53-KaPe-II}, Ovchinnikov \cite%
{54-Ov84}, Pelletier \cite{55-Pe} and the references listed there
(especially \cite{50-BrKrSe}).}

This paper is a first modest step in this direction. We speak of the
"elementary theory" in order to emphasize that we really do not prove
anything of significance. What we do is broadly speaking that we state a
program: we give a number of general definitions and point out some problems
motivated by these definitions. As is already apparent from the above we
rely heavily on the idea of category. In this connection the following
citation from Semadeni \cite{41-Se} might be useful to have in mind: "What
is significant in the theory of categories is the method of thinking, not
the particular results". (Besides this paper by Semadeni \cite{41-Se} let us
mention the following papers where categories are used in connection with
Banach spaces: Gelfand-Silov \cite{12-GeSi}$.$ Mitjagin-Svarc \cite{23-MiSv}$%
,$Semadeni \cite{41-Se} (almost identical with Semadeni \cite{42-Se} but
more detailed.). Also we have been inspired by parallel (structure) theory
of Banach spaces. It should of course be noted that the latter theory itself
is \underline{\emph{helas}}\emph{\ }far from complete. (We recommend the
reader in particular Lindenstrauss brilliant notes \cite{17-Li})$.\medskip $

Some ideas of this paper where briefly outlined in Peetre \cite{24-Pe}$.$
Most of the material was also treated in a seminar in Lund in december 1969. 
$\medskip $

The plan of the paper is as follows. In section 1 the notion of Banach
couple is introduced. In Section 2 we list some concrete Banach couples
which we need in the sequel as illustrations of the general development. In
Section $3$ we introduce several categories of Banach couples. In Section 4
we define subcouples and (the dual notion) quotient couples. Section 5 is
devoted to the important notion of retract. In Section 6 we discuss the
classification problem. In Section 7 we define the notions universal and
couniversal and prove an analogue of the Banach-Mazur theorem. In Section 8
we define the notions injective and projective. In particular the scalar
couple is seen as injective, in other words we have an analogue of the
Hahn-Banach theorem. We also introduce a general type, called
pseudo-injective. In Section 9 we apply our ideas to the problem of
determining all interpolation spaces. Here we have maybe the most
significant results of our investigation. Section 10 is devoted to a remark
on so-called pseudo-projective couples. Finally, Section 11 contains a
result concerning interpolation functions.

\section*{1. The notion of Banach couple}

A \underline{\emph{Banach couple}} consists of a Hausdorff topological
vector space $\mathcal{A}$ together with two Banach spaces $A_{0}$ and $%
A_{1} $ contained in $\mathcal{A}$ so that the corresponding embeddings $%
i_{0}$ and $i_{1}$ are linear continuous. (We denote the fields of scalars
by $\Lambda $, thus $\Lambda =\mathbb{R}$ or $\Lambda =\mathbb{C}$.) Thus we
have the diagram%
\begin{equation}
\begin{tabular}{ccc}
$A_{0}$ &  &  \\ 
& $\searrow $ &  \\ 
&  & $\mathcal{A}$ \\ 
& $\nearrow $ &  \\ 
$A_{1}$ &  & 
\end{tabular}
\tag{1.1}
\end{equation}%
We shall employ the notation $\overrightarrow{A}=\left\{ A_{0},A_{1}\right\} 
$ first used in Peetre \cite{35-Pe}. Thus neither the "large" space $%
\mathcal{A}$ nor the embeddings $i_{0}$ and $i_{1}$, will enter explicitly.
In most case one could assume that $\mathcal{A=}\Sigma \left( 
\overrightarrow{A}\right) $ where, generally speaking, $\Sigma \left( 
\overrightarrow{A}\right) =A_{0}+A_{1}=$ \underline{\emph{linear hull}}\emph{%
\ }of $A_{0}$ and $A_{1}$ but we shall not make this extra requirement. $%
\medskip $

Along with $\Sigma \left( \overrightarrow{A}\right) $ another subspace of $%
\mathcal{A}$ will be of particular interest, namely $\Delta \left( 
\overrightarrow{A}\right) =A_{0}\cap A_{1}=$\underline{\emph{intersection}}%
\emph{\ }of $A_{0}$ and $A_{1}.$We have $\Delta \left( \overrightarrow{A}%
\right) \subseteq A_{0}\subseteq \Sigma \left( \overrightarrow{A}\right) $
and $\Delta \left( \overrightarrow{A}\right) \subseteq A_{1}\subseteq \Sigma
\left( \overrightarrow{A}\right) .\medskip $

The notion of Banach couples can be generalized in several directions:

\subsection*{a.}

We could consider instead of just two spaces several spaces. We are then
lead to the notion of \underline{\emph{Banach }$\left( n+1\right) $-\emph{%
tuple}}\emph{, }where is any integer $\geq 0.$ For a Banach $\left(
n+1\right) -$ tuple we might use the notation $\overrightarrow{A}=\left\{
A_{0},A_{1},...,A_{n}\right\} .$ If $n=1$ we are back in the case of Banach
couples. In what follows we shall mostly be concerned with that special case
only. Sometimes, for formal reasons, the trivial case $n=0$ will be referred
to too; for most purposes a 1-tuple can be identified with a single Banach
space. For the general case we refer to Sparr \cite{44-Sp}.

\subsection*{b.}

We might use instead of Banach spaces other types of topological vector
space. We are then lead to the analogous notions of \underline{\emph{Hilbert
couple}}\emph{, \underline{\emph{quasi-Banach couple}}, \underline{\emph{%
locally convex couple}} }etc. We can also entirely drop the linear
structure. A particular instance of some interest is then the notion of 
\underline{\emph{metric couple}}\emph{\ }(cf. Peetre \cite{26-Pe},
Gustavsson \cite{15-Gu}). Also the "hybrid" \underline{\emph{normed Abelian
couple}}\emph{\ }is worthwhile to take into consideration (cf. Peetre-Sparr 
\cite{38-PeSp}). Again in what follows we shall mostly concentrate on the
case of Banach couples.

\subsection*{c.}

It is further sometimes necessary to let $A_{0}$ and $A_{1}$ be independent
beings with no direct connection with $\mathcal{A}$ and $i_{0}$ and $i_{1}$
be arbitrary continuous linear mappings not even necessarily injective. We
shall coin he word \underline{\emph{weak Banach couple}}\emph{\ }for this
case. Thus a weak Banach couple $\overrightarrow{A}=\left\{
A_{0},A_{1}\right\} $ is nothing but a diagram of the form $\left(
1.1\right) .$ We can still define $\Sigma \left( \overrightarrow{A}\right) $
as a subspace of $\mathcal{A},$ namely as the linear hull of $i_{0}\left(
A_{0}\right) $ and $i_{1}\left( A_{1}\right) ,$ but $\Delta \left( 
\overrightarrow{A}\right) $ has to be defined as a subspace of $A_{0}\times
A_{1},$ namely the one which consists of those pairs $\left(
a_{0},a_{1}\right) \in A_{1}\times A_{1}$ such that $i_{0}\left(
a_{0}\right) =i_{1}\left( a_{1}\right) $ (cf. Gagliardo \cite{10-Ga}). Thus
we have also the diagram:%
\begin{equation*}
\begin{tabular}{ccccccc}
&  & $A_{0}$ &  &  &  &  \\ 
& $\nearrow $ &  & $\searrow i_{0}$ &  &  &  \\ 
$\Delta \left( \overrightarrow{A}\right) $ &  &  &  & $\Sigma \left( 
\overrightarrow{A}\right) $ & $\rightarrow $ & $\mathcal{A}$ \\ 
& $\searrow $ &  & $\nearrow i_{1}$ &  &  &  \\ 
&  & $A_{1}$ &  &  &  & 
\end{tabular}%
\end{equation*}%
Here we shall not study weak Banach couples systematically but it is often
very useful to translate things to the language of weak Banach couples and
the reader is urger to try to do this whenever possible.$\medskip $

Given any Banach couple $\overrightarrow{A}=\left\{ A_{0},A_{1}\right\} $ we
introduce a one-parameter family of equivalent norms in $\Sigma \left( 
\overrightarrow{A}\right) $ as follows (first appearing in Peetre \cite%
{37-Pe} (and in disguised form in Gagliardo \cite{10-Ga}!); cf.
Butzer-Berens \cite{6-BuBe67} and many other places):%
\begin{equation*}
K\left( t,a;\overrightarrow{A}\right) =\inf_{a=a_{0}+a_{1}}\left( \left\Vert
a_{0}\right\Vert _{A_{0}}+t\left\Vert a_{1}\right\Vert _{A_{1}}\right) \text{
if }a\in \Sigma \left( \overrightarrow{A}\right) 
\end{equation*}%
where $0<t<\infty .$ Similarly, we define a one-parameter family of
equivalent norms in $\Delta \left( \overrightarrow{A}\right) $ as follows:%
\begin{equation*}
J\left( t,a;\overrightarrow{A}\right) =\max \left( \left\Vert a\right\Vert
_{A_{0}},t\left\Vert a\right\Vert _{A_{1}}\right) \text{ if }a\in \Delta
\left( \overrightarrow{A}\right) .
\end{equation*}%
For technical reasons it is sometimes necessary to use a modified definition:%
\begin{eqnarray*}
K_{p}\left( t,a;\overrightarrow{A}\right)  &=&\inf_{a=a_{0}+a_{1}}\left(
\left\Vert a_{0}\right\Vert _{A_{0}}^{p}+t^{p}\left\Vert a_{1}\right\Vert
_{A_{1}}^{p}\right) ^{1/p}\text{ } \\
J_{p}\left( t,a;\overrightarrow{A}\right)  &=&\left( \left\Vert a\right\Vert
_{A_{0}}^{p}+t^{p}\left\Vert a\right\Vert _{A_{1}}^{p}\right) ^{1/p}
\end{eqnarray*}%
where $1\leq p\leq \infty .$ Clearly $K=K_{1},J=J_{\infty }.$(We adopt the
convention to drop those arguments which are not of importance in a given
context.). A comparison of $K_{p}$ for various $p$ is made in
Holmstedt-Peetre \cite{16-HoPe}. $\medskip $

By an \underline{\emph{intermediate space}}\emph{\ }corresponding to $%
\overrightarrow{A}$ we mean a Banach space $A$ such that $\Delta \left( 
\overrightarrow{A}\right) \subseteq A\subseteq \Sigma \left( \overrightarrow{%
A}\right) $ where the embeddings are linear continuous. Thus $A_{0}$ and $%
A_{1}$ are intermediate spaces and so are $\Delta \left( \overrightarrow{A}%
\right) $ and $\Sigma \left( \overrightarrow{A}\right) .$ (It is understood
that $\Delta \left( \overrightarrow{A}\right) $ has the topology given by
the equivalent norm $J_{p}\left( t,a;\overrightarrow{A}\right) $ and $\Sigma
\left( \overrightarrow{A}\right) $ the one given by the equivalent norm $%
K_{p}\left( t,a;\overrightarrow{A}\right) .$)

\section*{2. Examples of Banach couples.}

We now give a list of special Banach couples which will be used in what
follows

\subsection*{a. (Lebesque couples)}

Let $M$ be a locally compact space provided with a positive Radon measure $%
\mu .$For any positive $\mu $ measurable function $\omega ,$ any $p$ with $%
1\leq p\leq \infty $ and any Banach space $A,$ we denote with $L_{p}\left(
\omega ,a\right) $ the space of $\mu $ measurable functions $a:M\rightarrow A
$ such that 
\begin{equation*}
\left\Vert a\right\Vert _{L_{p}\left( \omega ,A\right) }=\left(
\tint\nolimits_{M}\left( \omega \left( m\right) \left\Vert a\left( m\right)
\right\Vert _{A}\right) ^{p}d\mu \left( m\right) \right) ^{1/p}<\infty 
\end{equation*}%
(with the usual interpretation if $p=\infty ).$ If $\omega =1$ or $A=\Lambda 
$ (= field of scalars), we omit the corresponding argument, writing thus $%
L_{p}\left( A\right) $ or $L_{p}\left( \omega \right) $ or just $L_{p},$ if
both these requirements are met. If $M$ is a discrete space and each point
carries mass $1$ we use $l$ in place $L;$ we obtain $l_{p}\left( \omega
,A\right) ,l_{p}\left( \omega \right) ,l_{p}\left( A\right) ,l_{p}.$ If $%
\overrightarrow{A}=\left\{ A_{0},A_{1}\right\} $ is any Banach couple we can
form the couples%
\begin{equation*}
\overrightarrow{L}_{\overrightarrow{p}}\left( \overrightarrow{\omega },%
\overrightarrow{A}\right) =\left\{ L_{p_{0}}\left( \omega _{0},A_{0}\right)
,L_{p_{1}}\left( \omega _{1},A_{1}\right) \right\} 
\end{equation*}%
where we use the notation $\overrightarrow{p}=\left( p_{0},p_{1}\right) ,%
\overrightarrow{\omega }=\left( \omega _{0},\omega _{1}\right) .$ In two
special cases we have a simple "explicit" formula for $K:$

$\left( i\right) :p_{0}=1,p_{1}=\infty ,\omega _{0}=\omega
_{1}=1,A_{0}=A_{1}=A$

$\left( ii\right) :p_{0}=p_{1}=p$

Namely holds (cf. Peetre \cite{28-Pe})%
\begin{equation}
K\left( t,a;\left\{ L_{1}\left( A\right) ,L_{\infty }\left( A\right)
\right\} \right) =\tint\nolimits_{0}^{t}a^{\ast }\left( s\right) ds 
\tag{2.1}
\end{equation}

where $\ast $ stands for "decreasing absolute rearrangement of",%
\begin{eqnarray}
&&K_{p}\left( t,a;\left\{ L_{p}\left( \omega _{0},A_{0}\right) ,L_{p}\left(
\omega _{1},A_{1}\right) \right\} \right)  \TCItag{2.2} \\
&=&\left( \tint\nolimits_{M}\left( \omega _{0}\left( m\right) K_{p}\left( 
\frac{\omega _{1}\left( m\right) t}{\omega _{0}\left( m\right) },a;%
\overrightarrow{A}\right) \right) ^{p}d\mu \left( m\right) \right)
^{1/p},\left( \frac{1}{p}+\frac{1}{p^{^{\prime }}}=1\right)  \notag
\end{eqnarray}

(In particular holds $\left( A_{0}=A_{1}=A\right) $%
\begin{eqnarray}
&&K_{p}\left( t,a;\left\{ L_{p}\left( \omega _{0},A\right) ,L_{p}\left(
\omega _{1},A\right) \right\} \right)   \TCItag{2.3} \\
&=&\left( \tint\nolimits_{M}\left( \left( \frac{1}{\omega _{0}^{p^{^{\prime
}}}}+\frac{1}{\left( t\omega _{1}\right) ^{p^{^{\prime }}}}\right)
^{-p/p^{^{\prime }}}\left\Vert a\left( m\right) \right\Vert _{A}\right)
^{p}d\mu \left( m\right) \right) ^{1/p}.  \notag
\end{eqnarray}

\subsection*{b. (Lorentz couples).}

Similar construction with Lorentz space $L_{pq}\left( \omega ,A\right) $
instead of Lebesque space $L_{p}\left( \omega ,A\right) .$

\subsection*{c. (differentiable couples).}

Let $M$ be an open subset of $\mathbb{R}^{d}$ or more generally a $d$
dimensional Riemannian manifold. We denote by $C^{0}$ be the space of
bounded continuous functions $a:M\rightarrow \Lambda $ and by $C^{k}$ the
space of those functions k'th derivatives are in $C^{0}.$ We may now form
the couples 
\begin{equation*}
C^{\overrightarrow{k}}=\left\{ C^{k_{0}},C^{k_{1}}\right\} 
\end{equation*}%
where we write $\overrightarrow{k}=\left\{ k_{0},k_{1}\right\} .$ If $%
\Lambda =R,M=\mathbb{R}^{d},k=0,k_{1}=1$ we have again an "explicit" formula
for $K:$Namely holds (cf. Peetre \cite{29-Pe})%
\begin{equation}
K_{\infty }\left( t,a;\left\{ C^{0},C^{1}\right\} \right) =\frac{1}{2}\Omega
^{\ast }\left( t,a\right) ,  \tag{2.4}
\end{equation}%
where now $^{\ast }$ stands for "least concave majorant of" and $\Omega $ is
the modulus of continuity, which is, in view of Holmstedt-Peetre \cite%
{16-HoPe}, equivalent to%
\begin{equation}
K\left( t,a:\left\{ C^{0},C^{1}\right\} \right) =\sup_{x,y\in \mathbb{R}^{d}}%
\frac{\left\vert a\left( x\right) -a\left( y\right) \right\vert }{2+\frac{%
\left\vert x-y\right\vert }{t}}  \tag{2.5}
\end{equation}%
(with $x=\left( x_{1},..,x_{d}\right) ,\left\vert x\right\vert =\left\vert
x_{1}\right\vert +...+\left\vert x_{d}\right\vert $).

\subsection*{d. (Lipschitz couples).}

Similar construction with Lipschitz space $Lip_{\alpha }$ instead of $C^{k}$

\subsection*{e. (Sobolev couple).}

The Sobolev space $W_{p}^{k}$ is obtained if we in the definition of $C^{k}$
substitute the maximum norm for the $L_{p}$ norm (with respect to the Haar
measure). We obtain the couples.%
\begin{equation*}
W_{\overrightarrow{p}}^{\overrightarrow{k}}=\left\{
W_{p_{0}}^{k_{0}},W_{p_{1}}^{k_{1}}\right\} ,
\end{equation*}%
with $\overrightarrow{k}=\left\{ k_{0},k_{1}\right\} ,\overrightarrow{p}%
=\left\{ p_{0},p_{1}\right\} .$

\subsection*{f. (Besov couples).}

If $M=\mathbb{R}^{d}$ (we consider that case only) the Besov space $%
B_{p}^{sq}$ can be defined as follows (cf. Peetre \cite{30-Pe},\cite{31-Pe}%
): Introduce the following sets in $\widehat{\mathbb{R}^{d}}$ (dual group):%
\begin{equation*}
U_{\nu }=\left\{ \xi \mid :2^{\nu -1}\leq \left\vert \xi \right\vert \leq
2^{\nu +1}\right\} ,\left( v=0,\overset{+}{-}1,\overset{+}{-}2,...\right) 
\end{equation*}%
and choose also a family ("partition") of rapidly decreasing (i.e. belonging
to $\mathcal{S)}$ functions $\phi _{v}$ $\left( v=0,\overset{+}{-}1,\overset{%
+}{-}2,...\right) $ with the following properties%
\begin{eqnarray*}
\text{support of }\mathcal{F\phi }_{\nu } &\subseteq &U_{\nu }. \\
\left\vert D^{\alpha }\mathcal{F\phi }_{\nu }\left( \xi \right) \right\vert 
&\leq &C_{\alpha }\left\vert \xi \right\vert ^{-\alpha }\text{ (}C_{\alpha }%
\text{ independent of }\nu \text{!).} \\
\dsum_{\nu =-\infty }^{\infty }\left\vert \mathcal{F\phi }_{\nu }\left( \xi
\right) \right\vert  &\geq &c>0,
\end{eqnarray*}%
$\mathcal{F}$ denoting the Fourier transform. Then we define $B_{p}^{sq}$ to
be the space of those tempered distributions (i.e. belonging to $\mathcal{S}%
^{^{\prime }})$ such that%
\begin{equation*}
\left\Vert a\right\Vert _{B_{p}^{sq}}=\left( \dsum_{\nu =-\infty }^{\infty
}\left( 2^{2\nu s}\left\Vert \phi _{v}\ast \alpha \right\Vert
_{L_{p}}\right) ^{q}\right) ^{1/q}<\infty ,
\end{equation*}%
$\ast $ denoting convolution. Note that $B_{\infty }^{\alpha \infty
}=Lip_{\alpha }:$ The corresponding couples are%
\begin{equation*}
B_{\overrightarrow{p}}^{\overrightarrow{s}.\overrightarrow{q}}=\left\{
B_{p_{0}}^{s_{0}q_{0}},B_{p_{1}}^{s_{1}q_{1}}\right\} 
\end{equation*}%
with $\overrightarrow{s}=\left\{ s_{0},s_{1}\right\} ,\overrightarrow{q}%
=\left\{ q_{0},q_{1}\right\} ,\overrightarrow{p}=\left\{ p_{0},p_{1}\right\}
.$

\subsection*{g. (Spectral couples).}

Let $E$ be any Banach space and $P$ a closed linear operator in $E$ with
domain $D\left( P\right) \subseteq E.$Then we may consider the couple $%
\left\{ E,D\left( P\right) \right\} .$ Some of the previous examples are
special case of such a couple:

$\left( i\right) :\left\{ L_{p}\left( \omega \,_{0},A\right) ,L_{p}\left(
\omega _{1},A\right) \right\} $ with $Pa\left( m\right) =\frac{\omega
_{1}\left( m\right) }{\omega _{0}\left( m\right) }a\left( m\right) $

$\left( ii\right) :C^{\overrightarrow{k}}$ if $d=1,$ with $P=D^{k}$ ($D$ $=$
derivative), $k=k_{1}-k_{0}$

$\left( iii\right) :W_{p}^{\overrightarrow{k}}=\left\{
W_{p}^{k_{0}},W_{p}^{k_{0}}\right\} $ if $d=1,$ with $P=D^{k}$

or if $1<p<\infty ,$ with $P=\left( \sqrt{-\Delta }\right) ^{k}$ ($\Delta $
= Laplacian)

$\left( iv\right) :B_{p}^{\overrightarrow{s}q}=\left\{
B_{p}^{s_{0}q},B_{p}^{s_{1}q}\right\} ,$with $P=\left( \sqrt{-\Delta }%
\right) ^{k}.$

\subsection*{h. (scalar couples).}

If $\Lambda $ is the field of scalars then $\overrightarrow{\Lambda }%
=\left\{ \Lambda ,\Lambda \right\} $ is a couple.

\subsection*{i. (operator couples, dual couples).}

Let $\overrightarrow{A},\overrightarrow{B}$ be two Banach couples, where we
assume $\overrightarrow{A}$ is \underline{\emph{regular}}\emph{\ }in the
sense that $\Delta \left( \overrightarrow{A}\right) $ is dense in both $A_{0}
$ and $A_{1}$ and hence $\Sigma \left( \overrightarrow{A}\right) .$Then we
can consider the couple%
\begin{equation*}
\overrightarrow{\mathcal{\overline{M}}}\left( \overrightarrow{A},%
\overrightarrow{B}\right) =\left\{ \overline{\mathcal{M}}\left(
A_{0},B_{0}\right) ,\overline{\mathcal{M}}\left( A_{1},B_{1}\right) \right\}
,
\end{equation*}%
where $\overline{\mathcal{M}}\left( A,B\right) $ generally speaking denotes
the set of continuous linear mappings $T:A\rightarrow B$ (cf. Section 3 for
this notation). This is indeed a couple in view of the regularity of $%
\overrightarrow{A}.$ $\overline{\mathcal{M}}\left( A_{0},B_{0}\right) $ can
be identified with a subspace of $\overline{\mathcal{M}}\left( \Delta \left( 
\overrightarrow{A}\right) ,\Sigma \left( \overrightarrow{B}\right) \right) $
and the same for $\overline{\mathcal{M}}\left( A_{1},B_{1}\right) .$ If we
specialize: $\overrightarrow{B}=\overrightarrow{\Lambda }=\left\{ \Lambda
,\Lambda \right\} ,$ we obtain the couple%
\begin{equation*}
\overrightarrow{A}^{^{\prime }}=\left\{ A_{0}^{^{\prime }},A_{1}^{^{\prime
}}\right\} 
\end{equation*}%
where $A^{^{\prime }}$ generally speaking denotes the dual of $A:A^{^{\prime
}}=\overline{\mathcal{M}}\left( A,\Lambda \right) .$ If we drop the
assumption of regularity on $A,$ then $\overrightarrow{\overline{\mathcal{M}}%
}\left( \overrightarrow{A},\overrightarrow{B}\right) $ respectively $%
\overrightarrow{A}^{^{\prime }}$ are only weak couples in the sense of
Sub-Section 1c. (Concerning interpolation of operator couples see Peetre 
\cite{33-Pe}, \cite{34-Pe}.)

\section*{3. Some categories of Banach couples.}

We take for granted that the reader is acquainted with the basic ideas
connected with categories, functors, natural transformations etc., so called
"general nonsense" in the sense of Lang (cf. e.g. Mc Lane \cite{20-Mc}).
Roughly speaking a category $\mathcal{C}$ arises if we take a certain class $%
\mathcal{S}$ of "spaces" and adjoin to this class a suitable class $\mathcal{%
M}$ of "mappings" between these "spaces", a "space" being identified with
its identity mapping, in such a manner that certain simple axioms are
fulfilled. In particular we require%
\begin{equation}
T\in \mathcal{C},S\in \mathcal{C}\Longrightarrow TS\in \mathcal{C}:X\in 
\mathcal{C}\Longrightarrow id_{X}\in \mathcal{C}  \tag{3.1}
\end{equation}%
Given any category $\mathcal{C}$ we say that $T:X\rightarrow Y$ belonging to 
$\mathcal{C}$ is an \underline{\emph{isomorphism}}\emph{\ }(within $\mathcal{%
C)}$ if the inverse mapping $T^{-1}:Y\rightarrow X$ exists and also belongs
to $\mathcal{C}$, and that $X,Y$ belonging to $\mathcal{C}$ are \underline{%
\emph{isomorphic}}\emph{\ }or \underline{\emph{equivalent}}\emph{\ }(within $%
\mathcal{C)}$ if there exists at least one isomorphism $T:X\rightarrow
Y.\medskip $

Now we shall construct various categories of Banach couples. First however
we have to define what is meant by a mapping from a Banach couple into
another. If $\overrightarrow{A}=\left\{ A_{0},A_{1}\right\} $ and $%
\overrightarrow{B}=\left\{ B_{0},B_{1}\right\} $ are Banach couples we say
that $T$ is a mapping from $\overrightarrow{A}$ into $\overrightarrow{B}$
(in brief: $\overrightarrow{A}\rightarrow \overrightarrow{B})$ if $T$ is a
mapping from $\Sigma \left( \overrightarrow{A}\right) $ into $\Sigma \left( 
\overrightarrow{B}\right) $ such that $T\left( A_{0}\right) \subseteq
B_{0},T\left( A_{1}\right) \subseteq B_{1}$. We have thus the following
commutative diagram:%
\begin{equation*}
\begin{tabular}{cccccc}
$A_{0}$ &  & $\overset{T\mid _{A_{0}}}{\rightarrow }$ & $B_{0}$ &  &  \\ 
& $\searrow $ &  &  & $\searrow $ &  \\ 
&  & $\Sigma \left( \overrightarrow{A}\right) $ & $\overset{T}{\rightarrow }$
&  & $\Sigma \left( \overrightarrow{B}\right) $ \\ 
& $\nearrow $ &  &  & $\nearrow $ &  \\ 
$A_{1}$ &  & $\overset{T\mid _{A_{1}}}{\rightarrow }$ & $B_{1}$ &  & 
\end{tabular}%
\end{equation*}%
since that it follows that $T\left( \Delta \left( \overrightarrow{A}\right)
\right) \subseteq \Delta \left( \overrightarrow{B}\right) ,$and we can
complete the diagram accordingly. In all cases, $\mathcal{C}$ is the class
of all couples and $\mathcal{M}$ the specified class of mappings. Then we
can list some categories of interest. (The reader is advised to think out
what all this corresponds to in the case of Banach spaces, i.e. 1-tuples.)$%
\medskip $

$\overline{\mathcal{C}}_{b}:$ The category of bounded mappings. We say that $%
T:\overrightarrow{A}\rightarrow \overrightarrow{B}$ is bounded if $K\left(
t,Ta;\overrightarrow{B}\right) \leq CK\left( t,a;\overrightarrow{A}\right) $
whenever $a\in \Sigma \left( \overrightarrow{A}\right) ,$some $C$ (depending
on $T$ ).$\medskip $

$\overline{\mathcal{C}}_{l}.$ The category of linear continuous mappings. We
say that $T:\overrightarrow{A}\rightarrow \overrightarrow{B}$ is linear
continuous if $T$ considered as mapping from $\Sigma \left( \overrightarrow{A%
}\right) $ into $\Sigma \left( \overrightarrow{B}\right) $ is linear and if $%
\left\Vert Ta\right\Vert _{B_{0}}\leq C\left\Vert a\right\Vert _{A_{0}}$%
whenever $a\in A_{0}$ and $\left\Vert Ta\right\Vert _{B_{1}}\leq C\left\Vert
a\right\Vert _{A_{1}}$ whenever $a\in A_{1},$for some $C$ (depending on $T$).%
$\medskip $

$\overline{\mathcal{C}}_{Lip}.$ The category of Lipschitz mappings. We say
that $T:\overrightarrow{A}\rightarrow \overrightarrow{B}$ is Lipschitzian if 
$K\left( t,Ta^{^{\prime }}-Ta^{^{\prime \prime }}\right) \leq CK\left(
t,a^{^{\prime }}-a^{^{\prime \prime }}\right) $ whenever $a^{^{\prime
}},a^{^{\prime \prime }}\in \Sigma \left( \overrightarrow{A}\right) ,$for
some $C$ (depending on $T$).

\underline{\emph{Remark 3.1}}. The above definitions are not affected if we
replace $K$ by $K_{p}$ $\left( 1\leq p\leq \infty \right) .$ Also the
constant $C$ becomes the same, in view of the results of Holmstedt-Peetre 
\cite{16-HoPe}.

It is clear that $T\in \overline{\mathcal{C}}_{l}\Longrightarrow T\in 
\overline{\mathcal{C}}_{b}$ and $T\in \overline{\mathcal{C}}_{Lip}.$Also it
is plain, under the additional assumption $T0=0$ $T\in \overline{\mathcal{C}}%
_{Lip}\Longrightarrow T\in \overline{\mathcal{C}}_{b}.$We summarize:%
\begin{eqnarray*}
\overline{\mathcal{C}}_{l} &\subseteq &\overline{\mathcal{C}}_{b} \\
\overline{\mathcal{C}}_{l} &\subseteq &\overline{\mathcal{C}}_{Lip} \\
\overline{\mathcal{C}}_{Lip} &\subseteq &\overline{\mathcal{C}}_{b},T0=0
\end{eqnarray*}%
$\medskip $Let $\overline{\mathcal{M}}_{b}\left( \overrightarrow{A},%
\overrightarrow{B}\right) ,\overline{\mathcal{M}}_{l}\left( \overrightarrow{A%
},\overrightarrow{B}\right) ,\overline{\mathcal{M}}_{Lip}\left( 
\overrightarrow{A},\overrightarrow{B}\right) $ be the set of mappings $T:%
\overrightarrow{A}\rightarrow \overrightarrow{B}$ belonging to $\overline{%
\mathcal{C}}_{b},\overline{\mathcal{C}}_{l},\overline{\mathcal{C}}_{Lip}.$
If we set 
\begin{equation*}
\left\Vert T\right\Vert =\left\Vert T\right\Vert _{\overline{\mathcal{M}}%
_{x}\left( \overrightarrow{A},\overrightarrow{B}\right) }=\inf C\text{ }%
\left( x=b,l,Lip\right) 
\end{equation*}%
where $C$ runs through all constants $C$ appearing in the corresponding
defining conditions, we get a norm. It can be shown that $\overline{\mathcal{%
M}}_{x}\left( \overrightarrow{A},\overrightarrow{B}\right) $ $\left(
x=b,l,Lip\right) $ in this way becomes a Banach space. From the definition
of category (cf. $\left( 3.1\right) )$ follows that each class of spaces $%
\overline{\mathcal{M}}_{x}\left( \overrightarrow{A},\overrightarrow{B}%
\right) $ $\left( x=b,l,Lip\right) $ has the following \underline{\emph{%
multiplicative property}}\emph{:}%
\begin{equation*}
T\in \overline{\mathcal{M}}_{x}\left( \overrightarrow{A},\overrightarrow{B}%
\right) ,S\in \overline{\mathcal{M}}_{x}\left( \overrightarrow{B},%
\overrightarrow{C}\right) \Longrightarrow ST\in \overline{\mathcal{M}}%
_{x}\left( \overrightarrow{A},\overrightarrow{C}\right) 
\end{equation*}%
$\medskip $But we have moreover a corresponding \underline{\emph{%
multiplicative norm inequality}}\emph{: }%
\begin{equation*}
\left\Vert ST\right\Vert _{\overline{\mathcal{M}}_{x}\left( \overrightarrow{A%
},\overrightarrow{C}\right) }\leq \left\Vert S\right\Vert _{\overline{%
\mathcal{M}}_{x}\left( \overrightarrow{B},\overrightarrow{C}\right)
}\left\Vert T\right\Vert _{\overline{\mathcal{M}}_{x}\left( \overrightarrow{A%
},\overrightarrow{B}\right) }
\end{equation*}%
$\medskip $Thus we have here an instance of what might be called a 
\underline{\emph{normed category}}\emph{. }By this we mean a category $%
\mathcal{C}$ for which the set $\mathcal{M}\left( A,B\right) $ of all
mappings $T:X\rightarrow Y$ for all $X,Y$ in $\mathcal{C}$ is a Banach
space, the class of all sets $\mathcal{M}\left( X,Y\right) $, $X,Y$ running
though all of $\mathcal{C}$ satisfying an analogous multiplicative norm
inequality. The theory of normed categories ought to be developed for its
own sake. Here we make only the following

\underline{\emph{Remark 3.2}}. In every normed category $\mathcal{C}$ we can
define the concept of \underline{\emph{normed left ideal}}\emph{\ }$\mathcal{%
I}$. By this we mean a class of Banach spaces $\mathcal{I}\left( X,Y\right) $%
, labelled by $X,Y$ running through all of $\mathcal{C}$, such that%
\begin{eqnarray*}
T &\in &\mathcal{I}\left( X,Y\right) ,S\in \mathcal{M}\left( Y,Z\right)
\Longrightarrow ST\in \mathcal{I}\left( X,Z\right) \\
\left\Vert ST\right\Vert _{\mathcal{I}\left( X,Z\right) } &\leq &\left\Vert
S\right\Vert _{\mathcal{M}\left( Y,Z\right) }\left\Vert T\right\Vert _{%
\mathcal{I}\left( X,Y\right) }
\end{eqnarray*}

$\medskip $In the same way we can define \underline{\emph{normed right,
normed two-sided ideal}}\emph{. }An important question for ideals is the 
\underline{\emph{continuation problem}}\emph{: }To "continue" a normed ideal
defined in a subcategory of $\mathcal{C}$ to the whole of $\mathcal{C}.$ In
the case of Banach spaces normed ideal have been much studied (cf. in
particular Pietcsch \cite{39-Pi} where notably the continuation problem is
mentioned). Examples of normed ideal is provided by the class of all compact
operators, all nuclear operators, all absolutely summing operators etc. In
the case of Banach couples nothing is known.

$\medskip $\underline{\emph{Sub-remark}}\emph{. }Note that the concept of
normed ideal (left, right, two-sided) might be considered a sort of
continuous analogue of the concept of \underline{\emph{property}}\emph{\ }of
category. This this we mean a truth-false or better $0-1$ valued function $%
\iota $ defined in $\mathcal{C}$ such that e.g. the following requirements
are fulfilled%
\begin{eqnarray*}
\iota \left( S+t\right) &\leq &\max \left( \iota \left( S\right) ,\iota
\left( T\right) \right) \text{ (in the case of an "additive" category)} \\
\iota \left( ST\right) &\leq &\iota \left( S\right) \iota \left( T\right)
=\min \left( \iota \left( S\right) ,\iota \left( T\right) \right) \\
\iota \left( X\right) &=&\iota \left( id_{X}\right)
\end{eqnarray*}%
$\medskip $The continuation problem in a normed category has an analogue
here too. One can e.g. show that every property (in the usual sense) of the
class of spaces $\mathcal{S}$ of $\mathcal{C}$ can be extended the whole of $%
\mathcal{C}$ and that there exists a minimal extension.

$\medskip $\underline{\emph{Sub-example}}\emph{. }Consider the usual
category of Banach spaces (taking as mappings to linear continuous mappings
between two Banach spaces or $x=l$ in the above notation). Then if we
consider the property for a Banach space to be of finite dimension the
minimal extension is the property for a continuous linear mapping to be of
finite rank.

$\medskip $In what follows we shall however mostly try to work with the
corresponding \underline{\emph{exact}}\emph{\ }categories $\mathcal{C}_{l},%
\mathcal{C}_{b},\mathcal{C}_{Lip}$ rather than $\overline{\mathcal{C}}_{l},%
\overline{\mathcal{C}}_{b},\overline{\mathcal{C}}_{Lip}$ themselves. These
are obtained in similar way, by composing on the mappings $T$ the
restriction $C\leq 1$, or equivalently $\left\Vert T\right\Vert \leq 1.$In
analogy with the above we denote the set of such mappings by $\mathcal{M}%
_{b}\left( \overrightarrow{A},\overrightarrow{B}\right) $ etc. Thus $%
\mathcal{M}_{x}\left( \overrightarrow{A},\overrightarrow{B}\right) $ is
(closed) convex subset of $\mathcal{\overline{M}}_{x}\left( \overrightarrow{A%
},\overrightarrow{B}\right) .$ We have here in an instance of what might be
called a \underline{\emph{convex category}}\emph{.}

$\medskip $Also in what follows we shall almost entirely be concerned with
the cases $b$ and $l$ only.

$\medskip $We conclude this section by recalling the important notion of 
\underline{\emph{interpolation functor}}\emph{. }Roughly speak, this is
method how to get a space from a couple. To be precise, let $\mathcal{C=C}%
_{x}$, $x=b,l,Lip$ and let $\mathcal{C}^{0}$ be the corresponding category
of 1-tuples (i.e. Banach spaces). By an interpolation functor we mean a
functor $F:\mathcal{C\rightarrow C}^{0}$ such that for each couple $%
\overrightarrow{A}\in \mathcal{C}$ the space $F\left( \overrightarrow{A}%
\right) \in \mathcal{C}^{0}$ satisfies%
\begin{equation*}
\Delta \left( \overrightarrow{A}\right) \subseteq F\left( \overrightarrow{A}%
\right) \subseteq \Sigma \left( \overrightarrow{A}\right) 
\end{equation*}%
where the inclusions are continuous and moreover natural transformations. If 
$T:\overrightarrow{A}\rightarrow \overrightarrow{B}$ is any mapping in $%
\mathcal{C}$ it follows that the mapping $F\left( T\right) $ in $\mathcal{C}%
^{0}$ is just the restriction $T\mid _{F\left( \overrightarrow{A}\right) }.$%
By abuse of language, we agree to design restrictions of $T$ also by the
same letter. Thus we have the following \underline{\emph{interpolation
property}}\emph{: }%
\begin{equation*}
T:\overrightarrow{A}\rightarrow \overrightarrow{B}\Longrightarrow T:F\left( 
\overrightarrow{A}\right) \rightarrow F\left( \overrightarrow{B}\right) 
\end{equation*}%
which also may express by saying that $F\left( \overrightarrow{A}\right) $
and $F\left( \overrightarrow{B}\right) $ are \underline{\emph{interpolation
spaces}}\emph{\ }with respect to $\overrightarrow{A}$ and $\overrightarrow{B}
$ . $\medskip $Conversely all interpolation spaces can essentially be
obtained in this way, by virtue of an important result by Aronzajn-Gagliardo 
\cite{1-AG65} (formulated for $x=l$ only). Thus we make other convenient
abuse of language and refer to $F\left( \overrightarrow{A}\right) $ as an
interpolation space. Note that in the case of exact category (the
"un-barred" case) all norms which appear above are $\leq 1.$

$\medskip $\underline{\emph{Example 3.1}}\emph{. }Trivial examples of
interpolation spaces are $A_{0}$ and $A_{1}.$ Equally trivial but more
fruitful are $\Delta \left( \overrightarrow{A}\right) $ and $\Sigma \left( 
\overrightarrow{A}\right) .$ Indeed, in each of these two spaces we can
consider a whole family of equivalent norms, viz. $K\left( t,a;%
\overrightarrow{A}\right) $ and $J\left( t,a;\overrightarrow{A}\right) $ $%
\left( 0<t<\infty \right) ,$ and exploiting this circumstance we can
construct two very general types of non-trivial interpolation spaces, namely 
$K-$ and $J-$spaces. For an introduction to the theory of $K-$ and $J-$
spaces we refer to Butzer-Berens \cite{6-BuBe67} and Peetre \cite{28-Pe}.

\section*{4. Subcouples and quotient couples.}

Let $\overrightarrow{B}=\left\{ B_{0},B_{1}\right\} $ be a Banach couple,
i.e. $B_{0}$ and $B_{1}$ are two Banach spaces both being subspaces of a
"large" Hausdorff topological vector space $\mathcal{B}$ with continuous
embeddings. We obtain a \underline{\emph{sub-couple}}\emph{\ }$%
\overrightarrow{A}=\left\{ A_{0},A_{1}\right\} $ of $\overrightarrow{B}$ if
we take $A_{0}$ as any subspace of $B_{0}$ provided with a stronger norm (in
the sense that the embedding (injection) of $A_{0}$ into $B_{0}$ is
continuous) and as $A_{1}$ any subspace of $B_{1}$ provided with a stronger
norm, the large space being the same, $\mathcal{A=B}$. We space of a $l-$%
sub-couple if the norm agree, i.e.%
\begin{eqnarray}
\left\Vert a\right\Vert _{A_{0}} &=&\left\Vert a\right\Vert _{B_{0}}\text{
if }a\in A_{0}\subseteq B_{0}  \TCItag{4.1} \\
\left\Vert a\right\Vert _{A_{1}} &=&\left\Vert a\right\Vert _{B_{1}}~\text{%
if }a\in A_{1}\subseteq B_{1}  \notag
\end{eqnarray}

We speak of b-subcouple if in addition%
\begin{equation}
K\left( t,a;\overrightarrow{A}\right) =K\left( t,b;\overrightarrow{B}\right) 
\text{ if~}a\in A_{0}+A_{1}\subseteq B_{0}+B_{1}  \tag{4.2}
\end{equation}%
Note that $\left( 4.1\right) $ always implies half of $\left( 4.2\right) ,$%
viz. the inequality%
\begin{equation}
K\left( t,a;\overrightarrow{A}\right) \geq K\left( t,a;\overrightarrow{B}%
\right) \text{ if }a\in A_{0}+A_{1}\subseteq B_{0}+B_{1}  \tag{4.3}
\end{equation}%
\underline{\emph{Example 4.1}}$.$ $\left\{ C^{0},C^{1}\right\} $ is $l-$%
equivalent (l-isomorphic) to a b-subcouple of $\left\{ l^{\infty }\left( 
\frac{1}{2}\right) ,l^{\infty }\left( \frac{1}{\left\vert x-y\right\vert }%
\right) \right\} .$

This follows at once if we compare the formulas $\left( 2.3\right) $ and $%
\left( 2.4\right) $ (the special case $E=\mathbb{\mathbb{R}},$ $p=\infty ,$ $%
\omega _{0}\left( x,y\right) =\frac{1}{2},$ $\omega _{1}\left( x,y\right) =%
\frac{1}{\left\vert x-y\right\vert },$ $M=\mathbb{R\times \mathbb{R)}}$ .
Indeed the mapping%
\begin{equation*}
T:a\left( x\right) \rightarrow b\left( x,y\right) =a\left( x\right) -a\left(
y\right) 
\end{equation*}%
defines apparently an $l$ isomorphism from $\left\{ C^{0},C^{1}\right\} $
onto a b-subcouple of $\left\{ l^{\infty }\left( \frac{1}{2}\right)
,l^{\infty }\left( \frac{1}{\left\vert x-y\right\vert }\right) \right\}
.\medskip $

The following almost trivial result shows why b-subcouples are of interest.

\underline{\emph{Proposition 4.1}}. Let $\overrightarrow{A}$ be a
b-subcouple of $\overrightarrow{B}$ and let $F$ be any K-interpolation
functor. Then with equality of norm%
\begin{equation}
F\left( \overrightarrow{A}\right) =F\left( \overrightarrow{B}\right) \cap
\Sigma \left( \overrightarrow{A}\right)  \tag{4.4}
\end{equation}

\underline{\emph{Proof}}\emph{. }Let $\Phi $ be the function norm
corresponding to $F$ (cf. Peetre \cite{28-Pe}). Then, by definition%
\begin{equation*}
a\in F\left( \overrightarrow{A}\right) \iff a\in \Sigma \left( 
\overrightarrow{A}\right) \text{ and }\Phi \left[ K\left( t,a;%
\overrightarrow{A}\right) \right] <\infty 
\end{equation*}%
and similarly for $\overrightarrow{B}.$ Comparing this with $\left(
4.2\right) $ we readily get $\left( 4.4\right) .$\# $\medskip $

For l-subcouples we have only the inclusion%
\begin{equation}
F\left( \overrightarrow{A}\right) \subseteq F\left( \overrightarrow{B}%
\right) \cap \Sigma \left( \overrightarrow{A}\right)   \tag{4.5}
\end{equation}%
which follows in the same way from $\left( 4.3\right) .$ It is instructive
to look at the following counter-example.

\underline{\emph{Example 4.2}}. (Arne Persson). Take as $\overrightarrow{B}$
the couple $\left\{ L_{1},L_{\infty }\right\} $ with $M=T=$ torus (circle
group), $\mu =$ Haar measure. Choose as $\overrightarrow{A}$ the couple
obtained by restricting oneself to those functions $a$ on $T^{1}$ which have 
\underline{\emph{lacunary Fourier series}}\emph{, }i.e. 
\begin{equation*}
a\left( x\right) =\dsum_{k=1}^{\infty }a_{k}e^{i\left( 2^{k}x\right) }
\end{equation*}%
Consider the sequence $\alpha =\left( a_{k}\right) _{k=1}^{\infty }.$ It is
known that (see Zugmund \cite{48-Zy} vol. 1. p. 203 and p. 247 respectively):

$1^{0}$ the norms $\left\Vert \alpha \right\Vert _{l_{2}}$ and $\left\Vert
a\right\Vert _{L_{p}}$, $1\leq p<\infty $ are equivalent.

$2^{0}$ the norms $\left\Vert \alpha \right\Vert _{l_{1}}$ and $\left\Vert
a\right\Vert _{L_{\infty }}$ are equivalent.

It is well-known that $F\left( \overrightarrow{A}\right) =F\left( \left\{
L_{1},L_{\infty }\right\} \right) =L_{p}$ if $\frac{1}{p}=\frac{1-\theta }{1}%
+\frac{\theta }{\infty }$ with $F$ defined as $F\left( \overrightarrow{C}%
\right) =\overrightarrow{C}_{\theta p;K}$ for a general couple $%
\overrightarrow{C}.$ On the other hand since $F\left( \left\{
l_{2},l_{\infty }\right\} \right) =l_{qp}$ if $\frac{1}{q}=\frac{1-\theta }{2%
}+\frac{\theta }{\infty }$ we see that $F\left( \overrightarrow{B}\right) $
consists of those functions having a lacunary Fourier expansion such that $%
\alpha \in l_{qp.}$. Since $l_{qp}\neq l_{2}$ we conclude that $\left(
4.4\right) $ is violated. $\medskip $

Now we come to the dual notion. Let $\overrightarrow{B}$ as above. We obtain
a \underline{\emph{quotient couple}}\emph{\ }$\overrightarrow{A}=\left\{
A_{0},A_{1}\right\} $ of $\overrightarrow{B}$ if we take as $A_{0}$ any
quotient space, say $B_{0}/C_{0}$ of $B_{0}$ (provided with a norm such that
the projection of $A_{0}$ onto $B_{0}$ is continuous) and as $A_{1}$ any
quotient, say $B_{1}/C_{1}$ of $B_{1},$the large space being $\mathcal{B}%
/\Delta \left( \overrightarrow{C}\right) $ where $\overrightarrow{C}=\left\{
C_{0},C_{1}\right\} ,$if we also assume that $C_{1}\cap B_{0}\subseteq
C_{0},C_{0}\cap B_{1}\subseteq C_{1}$ otherwise we will only have a weak
couple in the sense of Section 1c. We speak of a l-quotient couples if the
norms agree, i.e.%
\begin{eqnarray}
\left\Vert a\right\Vert _{A_{0}} &=&\inf_{c\in C_{0}}\left\Vert
b+c\right\Vert _{B_{0}}~\text{if }a\in b+C_{0}\in A_{0}  \TCItag{4.6} \\
\left\Vert a\right\Vert _{A_{1}} &=&\inf_{c\in C_{1}}\left\Vert
b+c\right\Vert _{B_{1}}\text{ if }a\in b+C_{1}\in A_{1}  \notag
\end{eqnarray}%
We speak of a \underline{\emph{b-quotient couple}}\emph{\ }if in addition%
\begin{equation}
K\left( t,a;\overrightarrow{A}\right) =\inf_{x\in C_{0}+C_{1}}K\left( t,b+c;%
\overrightarrow{B}\right) \text{ if }a=b+C_{0}+C_{1}\in A_{0}+A_{1} 
\tag{4.7}
\end{equation}%
Clearly the notions of subcouple and quotient-couple are dual to each other,
as usual. In what follows we have however not been able to develop the dual
theory to the full extent. What the reason for this is, whether there is
something hidden or not, is to clear to me at the time being.

\section*{5. Retracts}

Given any category $\mathcal{C}$ and of spaces $X,Y$ in this category $%
\mathcal{C}$ we say that $X$ is \underline{\emph{retract}}\emph{\ }of $Y$
(in $\mathcal{C)}$ if there exists mappings $\alpha :X\rightarrow Y,\beta
:Y\rightarrow X$ in $\mathcal{C}$ such that $\beta \circ \alpha =id.$ In
other words we have then the commutative diagram%
\begin{equation*}
\begin{tabular}{ccc}
$X$ & $\overset{\alpha }{\rightarrow }$ & $Y$ \\ 
$\downarrow id$ & $\swarrow \beta $ &  \\ 
$X$ &  & 
\end{tabular}%
\end{equation*}%
(We also say that $\beta $ is a\underline{ \emph{retraction}}\emph{\ }of $%
\alpha .).$This is clearly a transitive relation. If $X$ is retract of $Y$
and $Y$ is retract of $X$ then $X$ is retract of $Z.$Note also that $X$ is
trivially a retract of $Y$ if $X$ is equivalent to $Y$ (take $\beta =\alpha
^{-1}$ if $\alpha $ is the isomorphism!). However $X$ retract of $Y,Y$
retract of $X$ does not in general imply $X$ equivalent to $Y.$ It is
customary to say that $X$ and $Y$ belong to the same \underline{\emph{%
retract class}}\emph{. }

Now return to the case of Banach couples. The commutative diagram in our
notation now takes the form:%
\begin{equation*}
\begin{tabular}{ccc}
$\overrightarrow{A}$ & $\overset{\alpha }{\rightarrow }$ & $\overrightarrow{B%
}$ \\ 
$\downarrow id$ & $\swarrow \beta $ &  \\ 
$\overrightarrow{A}$ &  & 
\end{tabular}%
\end{equation*}%
We have a different notion for each of our two categories $\mathcal{C}_{b}$
and $\mathcal{C}_{l}$ (the case of $\mathcal{C}_{Lip}$ is excluded here!).
However it turns out that one also has to make the two classes $b$ and $l$
interplay. Therefore we really have \underline{\emph{four}}\emph{\ }notions
of retract to keep in mind:

$\overrightarrow{A}$ is a \underline{\emph{b-retract}}\emph{\ }of $%
\overrightarrow{B}$ if $\alpha \in \mathcal{C}_{b},\beta \in \mathcal{C}_{b}$

$\overrightarrow{A}$ is a \underline{\emph{l-retract}}\emph{\ }of $%
\overrightarrow{B}$ if $\alpha \in \mathcal{C}_{l},\beta \in \mathcal{C}_{l}$

$\overrightarrow{A}$ is a \underline{\emph{lb-retract}}\emph{\ }of $%
\overrightarrow{B}$ if $\alpha \in \mathcal{C}_{l},\beta \in \mathcal{C}_{b}$

$\overrightarrow{A}$ is a \underline{\emph{bl-retract}}\emph{\ }of $%
\overrightarrow{B}$ if $\alpha \in \mathcal{C}_{b},\beta \in \mathcal{C}_{l}$

\underline{\emph{Example 5.1}}$.$ If $\overrightarrow{A}$ is $b-$subcouple
of $\overrightarrow{B}$, then $\overrightarrow{A}$ is a lb-retract of $%
\overrightarrow{B}.$ Namely if $\alpha =$ injection of $\overrightarrow{A}$
into $\overrightarrow{B},$ the corresponding retraction $\beta $ of $\alpha $
can be constructed as follows:%
\begin{eqnarray*}
\beta \left( b\right) &=&b\text{ if }b\in \Sigma \left( \overrightarrow{A}%
\right) \\
\beta \left( b\right) &=&0\text{ if }b\notin \Sigma \left( \overrightarrow{A}%
\right)
\end{eqnarray*}%
Thus we get a new way of looking at the definitions of Section 4.

\underline{\emph{Example 5.2}}. If $\overrightarrow{A}$ is a l-retract of $%
\overrightarrow{B}$, then as readily seen, $\overrightarrow{A}$ is
b-subcouple of $\overrightarrow{B}.$ Thus prop. 4.1 is applicable in this
case (cf. Baouendi-Goulaouic \cite{3-BaCo}, Baiocchi \cite{2-Ba66}).

\underline{\emph{Example 5.3}}. $B_{\overrightarrow{p}}^{\overrightarrow{s}%
\overrightarrow{q}}$ is a l-retract of $l_{\overrightarrow{q}}\left( 2^{%
\overrightarrow{s}v},L_{\overrightarrow{p}}\right) $ (the bared $\left( 
\overline{}\right) $ case!). Indeed take 
\begin{eqnarray*}
\alpha  &:&a\rightarrow \left( \phi _{v}\ast \alpha \right)  \\
\beta  &:&\left( a_{v}\right) \rightarrow \dsum_{v=-\infty }^{\infty }\psi
_{v}\ast \alpha _{v}
\end{eqnarray*}%
where $\phi _{v}$ is the partition from Sub-Section 2f and $\psi _{v}$ is
second partition chosen such that 
\begin{equation*}
\dsum_{v=-\infty }^{\infty }\mathcal{F\psi }_{v}\mathcal{F\phi }_{v}=1.
\end{equation*}

\underline{\emph{Remark 5.1}.} The same applies of course also for the space
(= 0-tuple) $B_{p}^{sq}.$Indeed, it is possible to prove a much sharper
result namely that $B_{p}^{sq}$ is isomorphic to $l_{q}\left( l_{p}\right) $
(Peetre \cite{30-Pe}).

\underline{\emph{Problem 5.1}\textbf{.}} To sharpen ex. 5.1 in the direction
or remark 5.1.

\underline{\emph{Problem 5.2}}\textbf{.} Is $\left\{ C_{0},C_{1}\right\} $ a
l-retract of $\left\{ l^{\infty }\left( \frac{1}{2}\right) ,l^{\infty
}\left( \frac{1}{\left\vert x-y\right\vert }\right) \right\} $ (bared or
unbared case does not matter)?. That have a $lb$ retract if plain if we
combine ex. 5.1 with ex. 4.1.

\underline{\emph{Problem 5.3}}. Same problem for $\left\{ Lip_{\alpha
_{0}},Lip_{\alpha _{1}}\right\} :\medskip $

For technical reasons it become also useful to introduce a hybrid which we
term \underline{\emph{pseudoretract}}. Roughly speaking it a case and $%
\alpha $ is in a way partially defined. Here is the exact corresponding to
one of our previous four case: We say that $\overrightarrow{A}$ is
bl-pseudoretract of $\overrightarrow{B}$ if for every $a\in \Sigma \left( 
\overrightarrow{A}\right) $ there exists $b\in \Sigma \left( \overrightarrow{%
B}\right) $ with $K\left( t,a;\overrightarrow{A}\right) =K\left( t,b;%
\overrightarrow{B}\right) $ and a mapping $\beta :\overrightarrow{B}%
\rightarrow \overrightarrow{A}\in \mathcal{C}_{1}$ (depending on $a!)$ such
that $b=\beta a.$ It is plain that every bl-retract is a bl-pseudoretract.
(Take $b=\alpha a$!) Let us hasten to give an example which indicates the
usefulness of the concept. $\medskip $

\underline{\emph{Example 5.4}}\emph{. }Let $\overrightarrow{A}$ be any
couple. We say that $a\in \Sigma \left( \overrightarrow{A}\right) $ is 
\underline{\emph{regular}}\emph{\ }if it belongs to the closure (in the
topology of $\Sigma \left( \overrightarrow{A}\right) )$ of $\Delta \left( 
\overrightarrow{A}\right) .$ One can show that $a\in \Sigma \left( 
\overrightarrow{A}\right) $ is regular iff 
\begin{equation}
K\left( t,a\right) =O\left( \max \left( 1,t\right) \right) ,t\rightarrow 0%
\text{ or }\infty .  \tag{$5.1$}
\end{equation}%
\emph{\ }According to the "fundamental lemma of the theory of interpolation
spaces" (cf. Peetre \cite{28-Pe}) there exists a constant $c$ (depending on $%
\overrightarrow{A})$ such that if $a\in \Sigma \left( \overrightarrow{A}%
\right) $ is regular and $\epsilon >0$ we have the integral representation%
\begin{equation}
a=\tint\nolimits_{0}^{\infty }u\left( t\right) \frac{dt}{t}  \tag{$5.2$}
\end{equation}%
for some (vectorvalued!) function $u\left( t\right) $ with%
\begin{equation}
J\left( t,u\left( t\right) \right) \leq \left( c+\epsilon \right) K\left(
t,a\right) .  \tag{$5.3$}
\end{equation}%
Let us set 
\begin{equation*}
\gamma \left( \overrightarrow{A}\right) =\inf c
\end{equation*}%
(which values $\gamma \left( \overrightarrow{A}\right) $ can attain is not
clear. I know only that $\inf_{\overrightarrow{A}}\gamma \left( 
\overrightarrow{A}\right) =\frac{1}{2}$ and $\sup_{\overrightarrow{A}}\gamma
\left( \overrightarrow{A}\right) <\infty .$Cf. Peetre \cite{32-Pe}. \cite%
{35-Pe}.) We claim that the functional $\gamma \left( \overrightarrow{A}%
\right) $ has the following \underline{\emph{\emph{monotonicity}}}\emph{%
\emph{\ }}$\overset{}{\emph{\emph{property}}}$\emph{\ }If $\overrightarrow{A}
$ is a bl-pseudoretract of $\overrightarrow{B}$ then 
\begin{equation}
\gamma \left( \overrightarrow{A}\right) \leq \gamma \left( \overrightarrow{B}%
\right)   \tag{5.4}
\end{equation}%
\emph{\ }Indeed assume that we have 
\begin{equation*}
b=\tint\nolimits_{0}^{\infty }v\left( t\right) \frac{dt}{t},J\left(
t,v\left( t\right) ;\overrightarrow{B}\right) \leq \left( c+\epsilon \right)
K\left( t,b;\overrightarrow{B}\right) 
\end{equation*}%
where $b\in \Sigma \left( \overrightarrow{B}\right) $ is the element which
corresponds to a given $a\in \Sigma \left( \overrightarrow{A}\right) $ in
the above definition of bl-pseudoretract. Note that if $a$ is regular so is $%
b$ (in view of $\left( 5.1\right) !).$ Now set 
\begin{equation*}
u\left( t\right) =\beta v\left( t\right) 
\end{equation*}%
where $\beta $ has the analogous position. Now follows readily%
\begin{eqnarray*}
\tint\nolimits_{0}^{\infty }u\left( t\right) \frac{dt}{t} &=&\beta \left(
\tint\nolimits_{0}^{\infty }v\left( t\right) \frac{dt}{t}\right) =\beta b=a
\\
J\left( t,u\left( t\right) ;\overrightarrow{A}\right)  &\leq &J\left(
t,v\left( t\right) ;\overrightarrow{B}\right) \leq \left( c+\epsilon \right)
K\left( t,b;\overrightarrow{B}\right) =\left( c+\epsilon \right) K\left( t,a;%
\overrightarrow{A}\right) 
\end{eqnarray*}%
i.e. $\left( 5.2\right) $ and $\left( 5.3\right) $ are implied. This clearly
establishes $\left( 5.4\right) .$

\section*{6. The classification problem.}

The classification problem poses itself for every "mathematical theory" (set
theory, group theory, theory of topological spaces etc.) once the notion of
equivalence (isomorphism) has been precise, i.e. when we have imbedded the
class of corresponding "spaces" $\mathcal{S}$ into a category $\mathcal{C}$.
It amounts essentially to the following:

a) Find all isomorphism classes.

A related problem, essentially contained in a), is

b) Find all retract classes.

$\medskip $(In some cases problem a) and b) are equivalent by a "Schr\"{o}%
der-Bernstein theorem".). $\medskip $

\underline{\emph{Example 6.1}}\emph{. }Let us consider problem a) and b) in
our category $\mathcal{C}_{b}$ (Section 3). Then both problems can be solved
if we can answer the following question: Which functions $\phi \left(
t\right) $ can appear as $K\left( t,a;\overrightarrow{A}\right) ?.$ To be
more specific, let us put, for any $\overrightarrow{A}$%
\begin{equation*}
S=S\left( \overrightarrow{A}\right) =\left\{ \phi :\exists a\in \Sigma
\left( \overrightarrow{A}\right) :K\left( t,a;\overrightarrow{A}\right)
=\phi \left( t\right) \right\} 
\end{equation*}%
$S$ is a cone: $\phi \in S,\lambda >0\Longrightarrow \lambda \phi \in S.$ It
is also clear that $S\subseteq S_{0}$ where $S_{0}$ denote the convex cone
of all positive concave functions $\phi .$ If $\overrightarrow{A}$ and $%
\overrightarrow{B}$ are b-isomorphic it follows that $S\left( 
\overrightarrow{A}\right) =S\left( \overrightarrow{B}\right) .$ Conversely
if $S\left( \overrightarrow{A}\right) =S\left( \overrightarrow{B}\right) $
and if moreover for any $\phi $ the inverse images of the mappings $%
a\rightarrow K\left( t,a;\overrightarrow{A}\right) $ and $b\rightarrow
K\left( t,b;\overrightarrow{B}\right) $ are equivalent, then $%
\overrightarrow{A}$ and $\overrightarrow{B}$ are b-isomorphic. A similar
statement for for the relation $\overrightarrow{A}$ is a b-retract of $%
\overrightarrow{B}.\medskip $

\underline{\emph{Sub-example}}\emph{. }Consider the couple $\left\{
L_{1},L_{\infty }\right\} .$Then $S=S\left( \left\{ L_{1},L_{\infty
}\right\} \right) $ is "almost" $S_{0},$the exact answer depending on $M.$%
The following cases are typical%
\begin{equation*}
\begin{tabular}{cc}
$M$ & $S$ \\ 
$\left( 0,\infty \right) $ & $\lim_{t\rightarrow 0}\phi \left( t\right) =0$
\\ 
$\left( 0,1\right) $ & $\lim_{t\rightarrow 0}\phi \left( t\right) =1,\phi
\left( t\right) =$ const. if $t\geq 1.$ \\ 
$\left\{ 0,1,2,3...\right\} $ & $\phi \left( t\right) $ is linear in each
subinterval $\left( i,i+1\right) $ $\left( i=0,1....\right) $%
\end{tabular}%
\end{equation*}

\underline{\emph{Problem 6.1}}. For which subsets $P\subseteq S_{0}$ does
there exists $\overrightarrow{A}\in \mathcal{C}$ such that $S\left( 
\overrightarrow{A}\right) =P.\medskip $

\underline{\emph{Remark 6.1}}\emph{. }$P$ need not to be closed in $S_{0}$%
which follows e.g. from the special case of the couple $\left\{
L_{p},L_{\infty }\right\} ,p>1.$\emph{\ }

\section*{7. Universal. Couniversal.}

A couple $\overrightarrow{B}$ is said to be \underline{\emph{b-universal}}
with respect to a subclass $\mathcal{S}_{1}$ of $\mathcal{S}$ (the class of
"spaces" of $\mathcal{S)}$ if for every $\overrightarrow{A}$ in $\mathcal{S}%
_{1},\overrightarrow{A}$ is b-equivalent to a b-subcouple of $%
\overrightarrow{B}.$It is plain that if $\overrightarrow{B}$ is b-universal
and $\overrightarrow{B}$ is b-subcouple of $\overrightarrow{C}$ then $%
\overrightarrow{C}$ is b-universal too.

\underline{\emph{Theorem 7.1}}. Every regular couple $\overrightarrow{A}$ is
equivalent to a b-subcouple of $\overrightarrow{l}_{\infty }=\left( 
\overrightarrow{\omega }\right) =\left\{ l_{\infty }\left( \omega
_{0}\right) ,l_{\infty }\left( \omega _{1}\right) \right\} $ for suitable $%
\omega _{0},\omega _{1},M.$

(It will follow from the proof below that is we put a restriction of the
cardinality of $\Sigma \left( \overrightarrow{A}\right) $ we can pick $%
\omega _{0},\omega _{1},M$ and thus $\overrightarrow{l}_{\infty }\left( 
\overrightarrow{\omega }\right) $ once for all. Thus $\overrightarrow{l}%
_{\infty }\left( \overrightarrow{\omega }\right) $ is then indeed
b-universal with respect to a suitable subclass of couples $\overrightarrow{A%
}.$

In other words, we have a "Banach-Mazur theorem".)

\underline{\emph{Proof}}: Consider the dual couple $\overrightarrow{%
A^{^{\prime }}}=\left\{ A_{0}^{^{\prime }},A_{1}^{^{\prime }}\right\} .$ We
take 
\begin{equation*}
\left\{ 
\begin{array}{c}
M=\text{ unit ball of }\Delta \left( \overrightarrow{A^{^{\prime }}}\right) 
\\ 
\omega _{0}\left( a^{^{\prime }}\right) =\left\Vert a^{^{\prime
}}\right\Vert _{A_{0}^{^{\prime }}} \\ 
\omega _{1}\left( a^{^{\prime }}\right) =\left\Vert a^{^{\prime
}}\right\Vert _{A_{1}^{^{\prime }}}%
\end{array}%
\right\} 
\end{equation*}%
and define a mapping $T:\overrightarrow{A}\rightarrow \overrightarrow{l}%
_{\infty }\left( \overrightarrow{\omega }\right) $:, $a\rightarrow f_{a}$
with $f_{a}\left( a^{^{\prime }}\right) =\left\langle a^{^{\prime
}},a\right\rangle $ where $\left\langle ,\right\rangle $ stands for the
duality. Now $K_{\infty }\left( t,a;\overrightarrow{A}\right) $ and $%
J_{1}\left( \frac{1}{t},a^{^{\prime }};\overrightarrow{A^{^{\prime }}}%
\right) $ are dual norms. It follows that (by the Hahn-Banach theorem and $%
\left( 2.2\right) )$%
\begin{eqnarray*}
K_{\infty }\left( t,a;\overrightarrow{A}\right)  &=&\sup \frac{\left\vert
\left\langle a^{^{\prime }},a\right\rangle \right\vert }{J_{1}\left( \frac{1%
}{t},a^{^{\prime }};\overrightarrow{A^{^{\prime }}}\right) }=\sup \frac{%
\left\vert f_{a}\left( a^{^{\prime }}\right) \right\vert }{\left\Vert
a^{^{\prime }}\right\Vert _{A_{0}^{-}}+\frac{1}{t}\left\Vert a^{^{\prime
}}\right\Vert _{A_{1}^{^{\prime }}}} \\
&=&K_{\infty }\left( t,f_{a};\overrightarrow{l}_{\infty }\left( 
\overrightarrow{\omega }\right) \right) 
\end{eqnarray*}%
Thus $T$ is an b-isomorphism of $\overrightarrow{A}$ onto a b-subcouple of $%
\overrightarrow{l}_{\infty }\left( \overrightarrow{\omega }\right) .$\#

The dual notion is \underline{\emph{b-couniversal}} is defined in a similar
way (involving quotientcouples, instead of subcouples).

\underline{\emph{Problem 7.1}}\emph{. }The analogue of th. 7.1 for
couniversal.

\section*{8. Injective. Projective.}

A couple $\overrightarrow{I}$ is said to be l-injective if for any two
couples $\overrightarrow{A}$ and $\overrightarrow{B}$ such that $%
\overrightarrow{A}$ is b-subcouple of $\overrightarrow{B}$ (in symbol $%
\overrightarrow{A}\overset{b}{\subseteq }\overrightarrow{B})$ and every
mapping $T:\overrightarrow{A}\rightarrow \overrightarrow{I}$ in $\mathcal{C}%
_{1}$ there exists a mapping $S:\overrightarrow{B}\rightarrow 
\overrightarrow{I}$ in $\mathcal{C}_{1}$ making the following diagram
commutative:%
\begin{equation*}
\begin{tabular}{ccc}
&  & $\overrightarrow{I}$ \\ 
& $\nearrow T$ & $\uparrow S$ \\ 
$\overrightarrow{A}$ & $\overset{b}{\subseteq }$ & $\overrightarrow{B}$%
\end{tabular}%
\end{equation*}%
If $\overrightarrow{I}_{1}$ is a l-retract of $\overrightarrow{I}$ and $%
\overrightarrow{I}$ is l-injective then $\overrightarrow{I}_{1}$ is
l-injective too.

\underline{\emph{Theorem 8.1}}\emph{. }Let $\overrightarrow{I}$ be a
l-injective couple such that $\overrightarrow{I}$ is equivalent to a
b-subcouple of $\overrightarrow{l}_{\infty }\left( \overrightarrow{\omega }%
\right) =\left\{ l_{\infty }\left( \omega _{0}\right) ,l_{\infty }\left(
\omega _{1}\right) \right\} $ for suitable $\omega _{0},\omega _{1},M.$ Then 
$\overrightarrow{I}$ is indeed a l-retract of \emph{\ }$\overrightarrow{l}%
_{\infty }\left( \overrightarrow{\omega }\right) .$ Conversely every
l-retract of $\overrightarrow{l}_{\infty }\left( \overrightarrow{\omega }%
\right) $ is l-injective.

\underline{\emph{Proof}}\emph{:} The direct part of the proof is trivial,
because if $\overrightarrow{I}$ is equivalent to a b-couple $\overrightarrow{%
A}$ of any couple $\overrightarrow{B}$ then $\overrightarrow{I}$ obviously
is a l-retract of $\overrightarrow{B}:$ Take $T:\overrightarrow{A}%
\rightarrow \overrightarrow{I}$ to an equivalence and extend it to $S:%
\overrightarrow{B}\rightarrow \overrightarrow{I}$. Let $i:\overrightarrow{A}%
\rightarrow \overrightarrow{B}$ the injection. Then we have the diagram%
\begin{equation*}
\begin{tabular}{ccc}
$\overrightarrow{I}$ &  &  \\ 
& $\searrow \alpha $ &  \\ 
$\downarrow $ &  & $\overrightarrow{B}$ \\ 
& $\swarrow \beta $ &  \\ 
$\overrightarrow{I}$ &  & 
\end{tabular}%
\end{equation*}%
with $\alpha =i\circ T^{-1},\beta =S.$ For the converse, it suffices, by the
proceeding remarks, to show that $\overrightarrow{l}_{\infty }\left( 
\overrightarrow{\omega }\right) $ is l-injective. One can also restrict to
the case when $M$ is singleton, i.e. $\overrightarrow{l}_{\infty }\left( 
\overrightarrow{\omega }\right) $ is the couple $\overrightarrow{\Lambda }_{%
\overrightarrow{\omega }}=\left\{ \Lambda _{\omega _{0}},\Lambda _{\omega
}\right\} $ where $\Lambda _{\omega }$ denotes $\Lambda $ with the norm $%
\left\Vert x\right\Vert =\frac{\left\vert x\right\vert }{\omega }.$Then $%
\overrightarrow{A}\rightarrow \overrightarrow{\Lambda }_{\overrightarrow{%
\omega }}$ is defined by a linear functional on $\Sigma \left( 
\overrightarrow{A}\right) $ such that 
\begin{equation*}
\left\vert Ta\right\vert \leq \omega _{0}\left\Vert a\right\Vert
_{A_{0}},\left\vert Ta\right\vert \leq \omega _{1}\left\Vert a\right\Vert
_{A_{1}}
\end{equation*}%
or%
\begin{equation*}
\left\vert Ta\right\vert \leq \omega _{0}K\left( \frac{\omega _{1}}{\omega
_{0}},a;\overrightarrow{A}\right) =\omega _{0}K\left( \frac{\omega _{1}}{%
\omega _{0}},a;\overrightarrow{B}\right) ;
\end{equation*}%
where the last equality follows from $\left( 4.2\right) .$By the Hahn-Banach
theorem follows that we can extend $T$ to a linear functional $S$ on $\Sigma
\left( \overrightarrow{B}\right) $ such that 
\begin{equation*}
\left\vert Sb\right\vert \leq \omega _{0}K\left( \frac{\omega _{1}}{\omega
_{0}},b;\overrightarrow{B}\right)
\end{equation*}%
But then follows also%
\begin{equation*}
\left\vert Sb\right\vert \leq \omega _{0}\left\Vert b\right\Vert
_{B_{0}},\left\vert Sb\right\vert \leq \omega _{1}\left\Vert b\right\Vert
_{B_{1}}
\end{equation*}%
so $S$ gives rise to a mapping $\overrightarrow{B}\rightarrow 
\overrightarrow{\Lambda }_{\overrightarrow{\omega }}.$ \#$\medskip $

The dual notion is \underline{\emph{l-projective}} and is defined in a
similar way (involving quotient couples instead of subcouples). We shall not
discuss it presently. Cf. however Section 10.$\medskip $

Here is however the definition of another hybrid., the introduction of which
is motivated by the above th. 8.1: We say that $\overrightarrow{I}$ is 
\underline{\emph{bl-pseudoinjective}} if $\overrightarrow{I}$ is a
bl-pseudoretract of some $\overrightarrow{l}_{\infty }\left( \overrightarrow{%
\omega }\right) .$

\section*{9. The problem of interpolation spaces.}

The general question we have in mind is the following: \underline{To
determine all interpolation spaces} or more precisely: $\medskip $

\underline{\emph{Problem 9.1}}\emph{. }Given any two couples $%
\overrightarrow{A}$ and $\overrightarrow{B}$ and spaces $A$ and $B$
satisfying%
\begin{equation*}
\Delta \left( \overrightarrow{A}\right) \subseteq A\subseteq \Sigma \left( 
\overrightarrow{A}\right) ,\Delta \left( \overrightarrow{A}\right) \subseteq
B\subseteq \Sigma \left( \overrightarrow{B}\right) 
\end{equation*}%
When does it hold true that%
\begin{equation*}
T:\overrightarrow{A}\rightarrow \overrightarrow{B}\in \mathcal{C}%
\Longrightarrow T:A\rightarrow B\in \mathcal{C}^{0}
\end{equation*}%
where $\mathcal{C=\overline{\mathcal{C}}}_{x}$ or $\mathcal{C=C}_{x}$ $%
x=b,l,Lip$ and $\mathcal{C}^{0}$ the corresponding category of 1-tuples
(Banach spaces), see Section 3. $\medskip $

By virtue of the theorem of Aronszajn-Gagliardo \cite{1-AG65} we see that
this is essentially equivalent to

\underline{\emph{Problem 9.2}}\emph{. }Let $\overrightarrow{A},%
\overrightarrow{B},\mathcal{C}$ problem 9.1 above and let $a\in \Sigma
\left( \overrightarrow{A}\right) ,b\in \Sigma \left( \overrightarrow{B}%
\right) .$When it is true that%
\begin{equation}
\exists T:\overrightarrow{A}\rightarrow \overrightarrow{B}\in \mathcal{C}%
:b=Ta  \tag{9.1}
\end{equation}

$\medskip $To fix the ideas consider the unbared case only. If $x=b,$i.e. $%
\mathcal{C=C}_{b},$ a necessary and sufficient condition for $\left(
9.1\right) $ to hold is%
\begin{equation}
K\left( t,b;\overrightarrow{B}\right) \leq K\left( t,a;\overrightarrow{A}%
\right)  \tag{9.2}
\end{equation}%
$\medskip $This is essentially the corresponding condition of Gagliardo \cite%
{11-Ga}. It is plain that if $x=l,$i.e. $\mathcal{C=C}_{l},$ then $\left(
9.2\right) $ is at least necessary too. Below we shall give a contribution
to the following

\underline{\emph{Problem 9.3}}\emph{. }Let $\overrightarrow{A},%
\overrightarrow{B}$ be as in problem 9.1, let $\mathcal{C=C}_{l},$let $a\in
\Sigma \left( \overrightarrow{A}\right) $ and $b\in \Sigma \left( 
\overrightarrow{B}\right) .$ When (= for which $\overrightarrow{A}$ and $%
\overrightarrow{B})$ does $\left( 9.2\right) $ imply $\left( 9.1\right)
.\medskip $

An essentially equivalent formulation is \underline{When are all
interpolation spaces $K$ spaces?} (There arises also the dual question: 
\underline{When are all interpolation spaces J-spaces}? of which question we
know very little.) $\medskip $

\underline{\emph{Example 9.1}}\emph{. }That the answer is yes in the special
case $\left\{ l_{1},l_{\infty }\right\} $ for a finite $M,$ is the effect of
a classical theorem of Hardy-Littlewood-Polya, on doubly stochastic matrices
(see e.g. \cite{21-Mi}). This result has been extended to the case $M$
denumerable and to the continuous case $\left\{ L_{1},L_{\infty }\right\} $
by many authors. In particular problem 9.2 gets a complete answer at least
if we restrict $A$ and $B$ somewhat (see Calderon \cite{7-Ca66}, Mitjagin 
\cite{22-Mi}) so in this case "almost all" interpolation spaces are K-spaces.%
$\medskip $

There is also a partial extension to the case $\left\{ l_{p},l_{\infty
}\right\} $ and $\left\{ L_{p},L_{\infty }\right\} $ (Cotlar \cite{8-Co},
Lorentz-Shimogaki \cite{19-LoSh}, cf. Bergh \cite{4-BE71}).$\medskip $

\underline{\emph{Remark 9.1}}. The Mitjagin-Calderon result was so-to-speak
the starting point of this investigation. Strangely this very result does
not fit well into the picture obtained below.$\medskip $

\underline{\emph{Example 9.2}}\emph{. }For not all couples $\overrightarrow{A%
}$ and $\overrightarrow{B}$ the K-spaces exhaust all interpolation cases the
simplest case is perhaps the case $\overrightarrow{A}=\overrightarrow{B}%
=\left\{ L_{p},W_{p}^{2}\right\} $ with $p\neq 2.$Indeed by complex
interpolation we see that $A=B=W_{p}^{1}$ is an interpolation space but
certainly not a K-space (use lacunary Fourier series, as in ex. 4.1!).$%
\medskip $

We have however the following positive result (a slightly more general form
of which is th. 9.2 below). \footnote{%
See also Cwikel-Peetre \cite[Theorem 4.1, Remark 4.1]{51-CwPe}, for more
recent developments.}

\underline{\emph{Theorem 9.1}}\emph{. }Let $\overrightarrow{B}$ be an
injective fulfilling the extra requirement of th. 8.1 and let $%
\overrightarrow{A}$ be regular. Then $\left( 9.2\right) $ implies $\left(
9.1\right) .\medskip $

We give two proofs:$\medskip $

\underline{\emph{First Proof}}: It is no restriction as is easily seen, to
take $\overrightarrow{B}=\overrightarrow{l}_{\infty }\left( \overrightarrow{%
\omega }\right) .$ By th. 7.1 we may assume that $\overrightarrow{A}$ is a
b-subcouple of $\overrightarrow{B}.$ For simplicity take $M=\left( 0,\infty
\right) ,\omega _{0}\left( x\right) =1,\omega _{1}\left( x\right) =1/x.$ Let 
$K\left( t,b\right) \leq K\left( t,a\right) .$Then we can find a mapping $%
T_{1}:\overrightarrow{B}\rightarrow \overrightarrow{B}$ in $\mathcal{C}_{l}$
such that $b=T_{1}a.$To this end we observe that by $\left( 2.9\right) $ and
the results of Holmstedt-Peetre \cite{16-HoPe} generally speaking we have $%
K\left( t,a\right) =\left\vert a\right\vert ^{\ast },$with $\ast $ as in
Sub-Section 2c, so that in our case $\left\vert b\right\vert \leq \left\vert
a\right\vert ^{\ast }$. It is easy to reduce to the case $b=\left\vert
a\right\vert ^{\ast }$ ,$a$ continuous $\geq 0.$Then we can just define $%
T_{1}\overline{a}$ (with $\overline{a}$ generic) by%
\begin{equation*}
T_{1}\overline{a}\left( x\right) =\left\{ 
\begin{array}{c}
\overline{a}\left( x\right) \text{ if }a\left( x\right) =a^{\ast }\left(
x\right)  \\ 
\overline{a}\left( x_{1}\right) +\frac{x-x_{1}}{x_{2}-x_{1}}\overline{a}%
\left( x_{2}\right) ~\text{if }\left( x_{1},x_{2}\right) \text{ is the
largest interval }I\text{ } \\ 
\text{containing }x\text{ such that }a\left( \xi \right) <a\left( x\right)
,\xi \in I%
\end{array}%
\right\} 
\end{equation*}%
Having defined $T_{1}$ we get a mapping $T:\overrightarrow{A}\rightarrow 
\overrightarrow{B}$ with $b=Ta$ simply by taking restrictions. \#$\medskip $

\underline{\emph{Second proof}}\emph{: }This proof is seemingly much longer.
We start with a general analysis which is useful also in other connections.
Let $\overrightarrow{A}$ be regular, as in th 9.1, but let $\overrightarrow{B%
}$ be a couple which is dual in the sense of Sub-Section 2i for some regular
couple $\overrightarrow{C}=\left\{ C_{0},C_{1}\right\} ,$i.e. $%
\overrightarrow{B}=\overrightarrow{C^{^{\prime }}}$. We denote duality by $%
\left\langle ,\right\rangle .$Let $b=Ta$ with $a\in \Sigma \left( 
\overrightarrow{A}\right) ,b\in \Sigma \left( \overrightarrow{B}\right) ,T:%
\overrightarrow{A}\rightarrow \overrightarrow{B}.$ If $c\in \Sigma \left( 
\overrightarrow{C}\right) $ we obtain%
\begin{equation*}
\left\langle b,c\right\rangle =\left\langle Ta,c\right\rangle =\left\langle
T,a\widehat{\otimes }c\right\rangle .
\end{equation*}%
We now claim that%
\begin{equation}
\left\vert \left\langle b,c\right\rangle \right\vert \leq K\left( 1,a\otimes
c;A_{0}\widehat{\otimes }C_{0},A_{1}\widehat{\otimes }C_{1}\right)  \tag{9.3}
\end{equation}

where $\widehat{\otimes }$ stands for the projective tensor product (see
Grothendieck \cite{14-Gr}) Indeed we have for any $a\in \Sigma \left( 
\overrightarrow{A}\right) $ and $c\in \Sigma \left( \overrightarrow{C}%
\right) $%
\begin{equation*}
\left\vert \left\langle T,a\otimes c\right\rangle \right\vert \leq
\left\Vert a\right\Vert _{A_{0}}\left\Vert c\right\Vert _{C_{0}}
\end{equation*}%
which implies%
\begin{equation*}
\left\vert \left\langle T,u\right\rangle \right\vert \leq \left\Vert
u\right\Vert _{A_{0}\widehat{\otimes }C_{0}}
\end{equation*}%
for any $u\in A_{0}\widehat{\otimes }C_{0}.$In the same way we get

\begin{equation*}
\left\vert \left\langle T,u\right\rangle \right\vert \leq \left\Vert
u\right\Vert _{A_{1}\widehat{\otimes }C_{1}}
\end{equation*}%
for any $u\in A_{1}\widehat{\otimes }C_{1}.$Returning to our original $a$
and $c,$we conclude if $a\otimes c=u_{0}+u_{1}$ then%
\begin{equation*}
\left\vert \left\langle b,c\right\rangle \right\vert \leq \left\Vert
u_{0}\right\Vert _{A_{0}\widehat{\otimes }C_{0}}+\left\Vert u_{1}\right\Vert
_{A_{1}\widehat{\otimes }C_{1}}
\end{equation*}%
which clearly implies $\left( 9.3\right) .$ The import point is now that
conversly if $\left( 9.3\right) $ holds true it follows from the Hahn-Banach
theorem the existence of $T$ in $\left( 9.1\right) .$Now we take $%
\overrightarrow{B}=\overrightarrow{l}_{\infty }\left( \overrightarrow{\omega 
}\right) .$Then we can take $\overrightarrow{C}=\overrightarrow{l_{1}}\left(
1/\overrightarrow{\omega }\right) .$But generally speaking $A\widehat{%
\otimes }l_{1}\left( \omega \right) =l_{1}\left( \omega ,A\right) $ (see
Grothendieck \cite{14-Gr}). This in this case $\left( 9.4\right) $ is
equivalent with%
\begin{equation*}
\left\vert \left\langle b,c\right\rangle \right\vert \leq K\left( 1,a\otimes
c;\overrightarrow{l_{1}}\left( 1/\overrightarrow{\omega },\overrightarrow{A}%
\right) \right) 
\end{equation*}%
or, in view of $\left( 2.2\right) ,$%
\begin{equation*}
\left\vert \left\langle b,c\right\rangle \right\vert \leq \dsum_{m\in M}%
\frac{1}{\omega _{0}\left( m\right) }K\left( \frac{\omega _{0}\left(
m\right) }{\omega _{1}\left( m\right) },a;\overrightarrow{A}\right)
\left\vert c\left( m\right) \right\vert 
\end{equation*}%
But in this case 
\begin{equation*}
\left\langle b,c\right\rangle =\dsum_{m\in M}b\left( m\right) c\left(
m\right) 
\end{equation*}%
Therefore we also get the equivalent condition%
\begin{equation*}
\left\vert b\left( m\right) \right\vert \leq \frac{1}{\omega _{0}\left(
m\right) }K\left( \frac{\omega _{0}\left( m\right) }{\omega _{1}\left(
m\right) },a;\overrightarrow{A}\right) .
\end{equation*}%
Finally we note (cf. first proof) that%
\begin{equation*}
K\left( t,b\right) =\phi ^{\ast }\left( t\right) 
\end{equation*}%
where%
\begin{equation*}
\phi \left( t\right) =\sup \left\{ \omega _{0}\left( m\right) \left\vert
b\left( m\right) \right\vert :\frac{\omega _{0}\left( m\right) }{\omega
_{1}\left( m\right) }=t\right\} 
\end{equation*}%
Thus $\left( 9.3\right) $ and $\left( 9.2\right) $ are equivalent in this
case. \#$\medskip $

It is now easy to find condition on $\overrightarrow{B}$ which is necessary
and sufficient for the validity of the conclusion of th. 9.1$\medskip $

\underline{$\emph{Theorem\ 9.2}$}\emph{. }$\left( 9.2\right) $ implies $%
\left( 9.1\right) $ for every regular $\overrightarrow{A}$ if and only if $%
\overrightarrow{B}$ is bl-pseudeinjective (see Section 8).

\underline{\emph{Proof}}: Let $\overrightarrow{B}$ be bl-pseudoinjective,
i.e. $\overrightarrow{B}$ is a pl-pseudoretract of some $\overrightarrow{I}=%
\overrightarrow{l}_{\infty }\left( \overrightarrow{\omega }\right) .$Let $%
\left( 9.2\right) $ be fulfilled. We can find $i\in \Sigma \left( 
\overrightarrow{I}\right) $ with $K\left( t,b\right) =K\left( t,i\right) $
and a mapping $\beta :\overrightarrow{I}\rightarrow \overrightarrow{B}$ such
that $b=\beta i.$Since $K\left( t,i\right) \leq K\left( t,a\right) ,$we can
also (by th. 9.1) find a mapping $S:\overrightarrow{A}\rightarrow 
\overrightarrow{I}$ such that $i=Sa.$Then $T=\beta \circ S$ apparently
fulfills the requirement for $\left( 9.1\right) .$This proves the "if" part.
The converse is trivial (take $\overrightarrow{A}=\overrightarrow{I}!).$\#$%
\medskip $

\underline{$\emph{Example\ 9.3}$}\emph{. }$\overrightarrow{L}_{%
\overrightarrow{p}\infty }=\left\{ L_{p_{0}\infty },L_{p_{1}\infty }\right\} 
$ (Lorentz couple, see Sub-Section 2b) is bl-pseudoinjective in $\overline{%
\mathcal{C}}.$ To see this let $b\in \Sigma \left( \overrightarrow{L}_{%
\overrightarrow{p}\infty }\right) .$With no loss of generality, we may
assume $b=b^{\ast },$taking for simplicity $M=\left( 0,\infty \right) .$ It
is readily seen that%
\begin{equation*}
K\left( t,b\right) \approx \sup \frac{b\left( x\right) }{%
x^{-1/p_{0}}+t^{-1}x^{-1/p_{1}}}=K\left( t,b;\overrightarrow{L}_{\infty
}\left( x^{-1/\overrightarrow{p}}\right) \right) 
\end{equation*}%
Thus we get the desired "pseudo-retraction" by considering the natural
injection $\overrightarrow{L}_{\infty }\left( x^{1/\overrightarrow{p}%
}\right) \rightarrow \overrightarrow{L}_{\overrightarrow{p}\infty }.$If we
apply th. 9.2 adapted to the bared case we obtain the following result: Let $%
\overrightarrow{A}$ be any regular couple and $\overrightarrow{B}=%
\overrightarrow{L}_{\overrightarrow{p}\infty }.$ Let $A$ and $B$ as in
problem 9.1. Then $T:\overrightarrow{A}\rightarrow \overrightarrow{B}$
implies $T:A\rightarrow B$ iff $K\left( t,b\right) \leq K\left( t,a\right) ,$
where $a\in \Sigma \left( \overrightarrow{A}\right) ,b\in \Sigma \left( 
\overrightarrow{B}\right) ,$implies $\left\Vert b\right\Vert _{B}\leq
C\left\Vert a\right\Vert _{A}$ for some constant $C\geq 1.$In particular if $%
\overrightarrow{A}=\overrightarrow{L}_{\overrightarrow{p}1}=\left\{
L_{p_{0}1},L_{p_{1}1}\right\} $ the latter condition reads:%
\begin{eqnarray*}
\sup \frac{b^{\ast }\left( x\right) }{x^{-1/p_{0}}+tx^{-1/p_{1}}} &\leq
&\tint\nolimits_{0}^{\infty }\min \left( x^{1/p_{0}},tx^{1/p_{1}}\right)
a^{\ast }\left( x\right) \frac{dx}{x} \\
&\Longrightarrow &\left\Vert b\right\Vert _{B}\leq C\left\Vert a\right\Vert
_{A}
\end{eqnarray*}%
or, specializing still more, if $p_{0}=1,p_{1}=\infty :$%
\begin{equation*}
tb^{\ast }\left( t\right) \leq \tint\nolimits_{0}^{t}a^{\ast }\left(
x\right) dx\Longrightarrow \left\Vert b\right\Vert _{B}\leq C\left\Vert
a\right\Vert _{A}
\end{equation*}%
All this is of course connected with the Marcinkiewicz interpolation theorem
(we refer to Calderon \cite{7-Ca66}. Boyd \cite{5-Bo67}, Semenov \cite{43-Se}%
, Zippin \cite{47-Zi}, Shimogaki \cite{45-Sh}). $\medskip $

We conclude this section by giving another application of condition $\left(
9.3\right) $ appearing in the (second) proof of th. 9.1. $\medskip $

In $\left( 9.3\right) $ take $\overrightarrow{A}=\overrightarrow{l}%
_{1}\left( \overrightarrow{\omega }\right) ,\overrightarrow{C}=%
\overrightarrow{l}_{1}\left( 1/\overrightarrow{\omega }\right) ,$thus $%
\overrightarrow{B}=\overrightarrow{l}_{\infty }\left( \overrightarrow{\omega 
}\right) .$Then we get%
\begin{equation*}
\left\vert \dsum_{m}b\left( m\right) c\left( m\right) \right\vert \leq
\dsum_{m}\dsum_{n}\min \left( \omega _{0}\left( m\right) /\omega _{0}\left(
n\right) ,\omega _{1}\left( m\right) /\omega _{1}\left( n\right) \right)
\left\vert a\left( m\right) \right\vert \left\vert c\left( n\right)
\right\vert 
\end{equation*}%
which clearly is equivalent to%
\begin{equation*}
\left\vert b\left( m\right) \right\vert \leq \dsum_{n}\min \left( \omega
_{0}\left( m\right) /\omega _{0}\left( n\right) ,\omega _{1}\left( m\right)
/\omega _{1}\left( n\right) \right) \left\vert a\left( n\right) \right\vert 
\end{equation*}%
Thus $\left( 9.4\right) $ is in this case a necessary and sufficient
condition for the validity of $\left( 9.1\right) .\medskip $

We can generalize this a little bit. Take again $\overrightarrow{A}=%
\overrightarrow{l}_{1}\left( \overrightarrow{\omega }\right) $ but let $%
\overrightarrow{C}$ and thus $\overrightarrow{B}$ be general. Thus we arrive
at the condition%
\begin{equation}
\left\vert \left\langle b,c\right\rangle \right\vert \leq C\dsum_{m}\omega
_{0}\left( m\right) K\left( \frac{\omega _{1}\left( m\right) }{\omega
_{0}\left( m\right) },c;\overrightarrow{C}\right) \left\vert a\left(
m\right) \right\vert   \tag{9.5}
\end{equation}%
Using duality, we can see that $\left( 9.5\right) $ is equivalent to the
following condition:%
\begin{eqnarray}
\exists \beta  &:&M\rightarrow \Delta \left( \overrightarrow{B}\right)
:b=\dsum_{m}a\left( m\right) \beta \left( m\right)   \TCItag{9.6} \\
\frac{1}{\omega _{0}\left( m\right) }J\left( \frac{\omega _{0}\left(
m\right) }{\omega _{1}\left( m\right) },\beta \left( m\right) \right)  &\leq
&1  \notag
\end{eqnarray}%
(This should be compared with the definition of J-spaces, cf. example 5.4.)

\section*{10. A remark on pseudoprojective couples.}

The possibility of defining projective couples was briefly referred to in
Section 5 but has not been exploited at all here. Let us however mention
that the couple $\overrightarrow{l}_{1}\left( \overrightarrow{\omega }%
\right) $ is projective. In this Section we make however some remarks in a
more general notion, namely pseudoprojective couples, which remarks are
intimately connected with the discussion in Section 9 (see also ex. 5.4).
More precisely we say that a Banach couple $\overrightarrow{A}$ is
bl-pseudoprojective, if for regular $\overrightarrow{A},$every $a\in \Sigma
\left( \overrightarrow{A}\right) $ and some $\overrightarrow{l}_{1}\left( 
\overrightarrow{\omega }\right) ,$we can find $\lambda \in l_{1}\left( 
\overrightarrow{\omega }\right) $ with $K\left( t,a;\overrightarrow{A}%
\right) =K\left( t,\lambda ;\overrightarrow{l}_{1}\left( \overrightarrow{%
\omega }\right) \right) $ and a mapping $\beta :\overrightarrow{l}_{1}\left( 
\overrightarrow{\omega }\right) \rightarrow \overrightarrow{A}$ in $\mathcal{%
C}_{l}$ such $a=\beta \left( \lambda \right) .$The structure of the elements 
$\mathcal{M}\left( \overrightarrow{l}_{1}\left( \overrightarrow{\omega }%
\right) ,\overrightarrow{A}\right) $ is easily determined. Thus we see that
the existence of $\beta $ is equivalent to the existence of a function $%
u:M\rightarrow \Delta \left( \overrightarrow{A}\right) $ such that have the
representation%
\begin{equation*}
a=\dsum_{m}\lambda \left( m\right) u\left( m\right) ;
\end{equation*}%
with%
\begin{equation*}
\dsum_{m}\min \left( \omega _{0}\left( m\right) ,t\omega _{1}\left( m\right)
\right) \left\vert \lambda \left( m\right) \right\vert =K\left( t,a;%
\overrightarrow{A}\right) 
\end{equation*}%
If we specialize: $M=\left( 0,\infty \right) ,\omega _{0}\left( x\right)
=1,\omega _{1}\left( x\right) =\frac{1}{x}$, we see that this is the
continuous analogue of the representation $\left( 5.2\right) $ with $\left(
5.3\right) $ replaced by%
\begin{equation*}
\tint\nolimits_{0}^{\infty }\min \left( 1,\frac{t}{x}\right) J\left(
x,u\left( x\right) \right) \frac{dx}{x}\leq K\left( t,a\right) .
\end{equation*}%
All this in the unbared case. In the bared case we are lead to the inequality%
\begin{equation*}
\tint\nolimits_{0}^{\infty }\min \left( 1,\frac{t}{x}\right) J\left(
x,u\left( x\right) \right) \frac{dx}{x}\leq \left( c^{^{\prime }}+\epsilon
\right) K\left( t,a\right) 
\end{equation*}%
where $c^{^{\prime }}$ is a constant depending on $\overrightarrow{A}.$%
Corresponding to $\gamma \left( A\right) $ (see ex. 5.4) we are lead to put%
\begin{equation*}
\gamma ^{^{\prime }}\left( \overrightarrow{A}\right) =\inf c^{^{\prime }}
\end{equation*}%
(Cf. Peetre \cite{32-Pe},\cite{35-Pe}.) The following examples will perhaps
illuminate the problem mentioned in connection with remark 9.1. $\medskip $

\underline{\emph{Example 10.1}}\emph{. } $\left\{ L_{1},L_{\infty }\right\} $
is a bl-pseudoretract of $\left\{ L_{1}\left( x\right) ,L_{1}\left( 1\right)
\right\} $ in $\mathcal{C}$ and thus bl-pseudoprojective. More generally $%
\overrightarrow{L}_{\overrightarrow{p},1}=\left\{
L_{p_{0}1},L_{p_{1}q}\right\} $(Lorentz couple) is a bl-pseudoretract of $%
\overrightarrow{L}_{1}\left( x^{1/\overrightarrow{p}}\right) =\left\{
L_{1}\left( x^{1/p_{0}}\right) ,L\left( x^{1/p_{1}}\right) \right\} .$This
is the dual result of the result in ex. 9.3 and the proof is quite parallel. 
$\medskip $

\underline{\emph{Example 10.2}}\emph{. }$\left\{ \mathcal{S}_{1},\mathcal{S}%
_{\infty }\right\} $ is a bl-pseudoretract of $\left\{ l_{1},l_{\infty
}\right\} $ and thus bl-pseudoprojective by the discrete version of ex. 10.1
($\mathcal{S}_{p}$ denote the so-called p-trace class of operators in
Hilbert space; see Gohberg-Krein \cite{13-GoKr}, see also Triebel \cite%
{46-Tr}, Peetre \cite{35-Pe}, for the connection with interpolation spaces).

\section*{11. An application to interpolation functions.}

In this Section we consider the problem of determining all interpolation
spaces (Section 9) in the special case of spectral couples (Sub-Section 2g).
Thus we taken in problem 9.1 $\overrightarrow{A}=\left\{ E,D\left( P\right)
\right\} ,\overrightarrow{B}=\left\{ F,D\left( Q\right) \right\} $ and also $%
A=D\left( h\left( P\right) \right) ,B=D\left( h\left( Q\right) \right) $
where $h$ is a given scalar valued function. If $A$ and $B$ are
interpolation spaces in the sense of problem 9.1 we say that $h$ is an 
\underline{\emph{interpolation function}} with respect to $P$ and $Q.$Here
are some papers dealing with interpolation function for special types of
operators: Lions \cite{18-Li}, Foias-Lions \cite{9-FoLi}, Peetre \cite{36-Pe}%
, \cite{37-Pe}, Schechter \cite{40-Sc}. Very often one take $P=Q$ (and $E=F)$
and this we do also in what follows. Thus there is given at the onset only
one spectral couple $\overrightarrow{A}=\left\{ E,D\left( P\right) \right\}
. $But now we assume that $\overrightarrow{A}$ is l-retract of the spectral
couple $\overrightarrow{A^{\ast }}=\left\{ E^{\ast },D\left( P^{\ast
}\right) \right\} .$More precisely we assume that we have a commutative
diagram%
\begin{equation*}
\begin{tabular}{ccc}
$\left\{ E,D\left( P\right) \right\} $ &  &  \\ 
& $\searrow \alpha $ &  \\ 
$\downarrow id$ &  & $\left\{ E^{\ast },D\left( P^{\ast }\right) \right\} $
\\ 
& $\swarrow \beta $ &  \\ 
$\left\{ E,D\left( P\right) \right\} $ &  & 
\end{tabular}%
\end{equation*}%
where $\alpha $ and $\beta $ satisfy, beside the basic relations, the
conditions%
\begin{equation}
\alpha P=P^{\ast }\alpha ,\beta P^{\ast }=P\beta  \tag{11.1}
\end{equation}%
Each of the conditions $\left( 11.1\right) $ clearly imply%
\begin{equation}
P=\beta P^{\ast }\alpha  \tag{11.2}
\end{equation}%
$\medskip $

\underline{\emph{Example 11.1}}\emph{. }From ex. 5.3 we know that $B_{p}^{%
\overrightarrow{s}q}=\left\{ B_{p}^{s_{0}q},B_{p}^{s_{1}q}\right\} $ is a
retract of $l_{q}\left( 2^{\overrightarrow{s}v},L_{p}\right) =\left\{
l_{q}\left( 2^{s_{0}v},L_{p}\right) ,l_{q}\left( 2^{s_{1v}},L_{p}\right)
\right\} .$By Sub-Section 2g (cases $\left( i\right) $ and $\left( iv\right)
)$ we also know that these are spectral couples, the first with $P=\left( 
\sqrt{-\Delta }\right) ^{s_{0}-s_{1}},$the second with $P^{\ast }=2^{v\left(
s_{1}-s_{0}\right) }.$But $\left( 11.1\right) $ nor $\left( 11.2\right) $ is
fulfilled. However this can be remedied upon if we replace the original
operator $\left( \sqrt{-\Delta }\right) ^{s_{1}-s_{0}}$by an operator of the
type 
\begin{equation*}
a\rightarrow \dsum_{v=-\infty }^{\infty }2^{v\left( s_{1}-s_{0}\right) }\rho
_{v}\ast a
\end{equation*}%
where $\rho _{v}$ is a third partition with suitable reproducing properties
with respect to the partition $\phi _{v}$ and $\psi _{v}$entering in the
definition of $\alpha $ and $\beta $ (see Sub-Section 2f). $\medskip $

\underline{\emph{Example 11.2}}\emph{. }Similar examples with $\left\{
L_{p},W_{p}\right\} $ and $\left\{ L_{p}\left( l_{2}\right) ,L_{p}\left(
l_{2}\left( 2^{v}\right) \right) \right\} ,$ with $1<p<\infty ,$because we
have to use the Paley-Littlewood theorem (cf. Peetre \cite{30-Pe}. \cite%
{31-Pe}). $\medskip $

Returning to the previous general situation, we can prove the following. $%
\medskip $

\underline{\emph{Theorem 11.1}}\emph{. }Every interpolation function with
respect to $P^{\ast }$ is an interpolation function with respect to $P.$

\underline{\emph{Proof:}} Let $T:\overrightarrow{A}\rightarrow 
\overrightarrow{A}$ be given. We define $T^{\ast }:\overrightarrow{A^{\ast }}%
\rightarrow \overrightarrow{A^{\ast }}$ by setting $T^{\ast }=\alpha \circ
t\circ \beta .$Then by $\beta \circ \alpha =id$ follows $T=\beta \circ
T^{\ast }\circ \alpha .$ If $h$ is an interpolation function with respect to 
$P^{\ast },$we conclude%
\begin{equation*}
T^{\ast }:D\left( h\left( P^{\ast }\right) \right) \rightarrow D\left(
h\left( P^{\ast }\right) \right)
\end{equation*}%
On the other hand we have%
\begin{eqnarray*}
\alpha &:&D\left( h\left( P\right) \right) \rightarrow D\left( h\left(
P^{\ast }\right) \right) \\
\beta &:&D\left( h\left( P^{\ast }\right) \right) \rightarrow D\left(
h\left( P\right) \right)
\end{eqnarray*}%
\emph{\ }

\bigskip because in view of $\left( 11.1\right) $ we have at least formally%
\begin{equation*}
\alpha h\left( P\right) =h\left( P^{\ast }\right) \alpha ,\beta h\left(
P^{\ast }\right) =P^{\ast }h\left( \beta \right) .
\end{equation*}%
Thus we get upon composition%
\begin{equation*}
T:h\left( P\right) \rightarrow h\left( P\right) 
\end{equation*}%
and $h$ is an interpolation function with respect to $P.$\#

\end{document}